\documentclass[3p, authoryear, times, final]{elsarticle}

        \usepackage{preamble}
        \makeatletter
        \def\ps@pprintTitle{}
        \makeatother

\journal{}

\begin{document}

\begin{frontmatter}

\title{Supervised Large Neighbourhood Search for MIPs}

            \author[1]{Charly Robinson La Rocca\corref{cor1}}
            \ead{charly.robinson.la.rocca@umontreal.ca}

            \author[2]{Jean-François Cordeau}
            \ead{jean-francois.cordeau@hec.ca}

            \author[1]{Emma Frejinger}
            \ead{emma.frejinger@umontreal.ca}

            \cortext[cor1]{Corresponding author}

            \affiliation[1]{organization={Department of Computer Science and Operations Research, Université de Montréal},
            addressline={2900, boul. Édouard-Montpetit}, 
            city={Montreal},
            postcode={H3T 1J4}, 
            country={Canada}}

            \affiliation[2]{organization={Department of Logistics and Operations Management, HEC Montréal},
            addressline={3000 ch. de la Côte-Sainte-Catherine}, 
            city={Montreal},
            postcode={H3T 2A7}, 
            country={Canada}}

\begin{abstract}
    Large Neighbourhood Search (LNS) is a powerful heuristic framework for solving Mixed-Integer Programming (MIP) problems. However, designing effective variable selection strategies in LNS remains challenging, especially for diverse sets of problems. In this paper, we propose an approach that integrates Machine Learning (ML) within the destroy operator of LNS for MIPs with a focus on minimal offline training. We implement a modular LNS matheuristic as a test bench to compare different LNS heuristics, including our ML-enhanced LNS.  Experimental results on the MIPLIB 2017 dataset demonstrate that the matheuristic can significantly improve the performance of state-of-the-art solvers like Gurobi and SCIP. We conduct analyses on noisy oracles to explore the impact of prediction accuracy on solution quality.  Additionally, we develop techniques to enhance the ML model through loss adjustments and sampling routines. Our findings suggest that while random LNS remains competitive, our Supervised LNS (SLNS) outperforms other baselines and helps set the foundation for future research on ML for LNS methods that are both efficient and general. 
\end{abstract}
    

\end{frontmatter}

\newcommand{\V}{\mathcal{V}}
\newcommand{\Vr}{\mathcal{V}'}
\renewcommand{\P}{\mathcal{P}}
\renewcommand{\Pr}{\mathcal{P}'}
\newcommand{\K}{\mathcal{K}}
\newcommand{\kt}{k'}
\newcommand{\ko}{k^*}
\newcommand{\x}{\mathbf{x}}
\newcommand{\ct}{c^T}
\newcommand{\xv}{x_v}
\newcommand{\xvk}{x_v^k}
\newcommand{\xlp}{\tilde{\mathbf{x}}}
\newcommand{\xvklp}{\tilde{x}_v^k}

\newcommand{\xvkt}{x_{vt}^k}
\newcommand{\kvt}{k_{vt}}
\newcommand{\kp}{\mathbf{k}_{v}}
\newcommand{\X}{\mathcal{X}_p}
\newcommand{\tp}{T_p}
\renewcommand{\r}{r}
\newcommand{\bnb}{B\&B}
\newcommand{\xc}{\mathbf{x}^*}
\newcommand{\xinc}{\mathbf{x}_{\text{inc}}}
\newcommand{\dmax}{d_{\text{max}}}

\newcommand{\cplex}{\text{CPLEX}}
\newcommand{\pnf}{\text{PNF} }
\newcommand{\runtime}{\text{runtime}}
\newcommand{\solve}{\text{solve}}
\newcommand{\probe}{\text{probe}}
\newcommand{\select}{\text{select}}
\newcommand{\score}{\text{score}}
\newcommand{\predict}{\text{predict}}

\newcommand{\bigx}{\mathcal{X}}
\newcommand{\bigy}{\mathcal{Y}}
\newcommand{\bigdstar}{\mathcal{D}^*}

\newcommand{\pred}[1]{\bar{#1}}
\newcommand{\inc}[1]{{#1}'}

\newcommand{\loss}{\mathcal{L}}
\newcommand{\risk}{\mathcal{R}}
\newcommand{\bigomega}{\Omega}

\newcommand{\fancyf}{\psi}
\newcommand{\bigpsi}{\Psi}
\newcommand{\agg}{\mathcal{A}}
\newcommand{\update}{\mathcal{U}}

\newcommand{\fcn}{MCFND}
\newcommand{\xijk}{x_{ij}^k}
\newcommand{\yij}{y_{ij}}
\newcommand{\ryij}{\tilde{y}_{ij}}
\newcommand{\yb}{\boldsymbol{y}}
\newcommand{\xb}{\boldsymbol{x}}

\newcommand{\ry}{\tilde{\yb}}
\newcommand{\yss}{\yb^{\text{SS}}}
\newcommand{\yls}{\yb^{\text{LS}}}

\newcommand{\yssij}{y_{ij}^{\text{SS}}}
\newcommand{\ybss}{\boldsymbol{y}^{\text{SS}}}

\newcommand{\byij}{\bar{y}_{ij}}
\newcommand{\by}{\bar{\yb}}
\newcommand{\yinc}{\by}

\newcommand{\yp}{\yb'}

\newcommand{\cijk}{c_{ij}^k}
\newcommand{\cij}{c_{ij}}

\newcommand{\fij}{f_{ij}}
\newcommand{\uij}{u_{ij}}
\newcommand{\Xij}{X_{ij}}
\newcommand{\dpij}{\delta^p_{ij}}
\newcommand{\Paths}{P^k}
\newcommand{\onepath}{p}
\newcommand{\zkp}{z^k_p}

\newcommand{\np}{N_i^{+}}
\newcommand{\nminus}{N_i^{-}}
\newcommand{\Pf}{\mathcal{P}}
\newcommand{\Plin}{\widetilde{\mathcal{P}}}
\newcommand{\pij}{\rho_{ij}}

\newcommand{\binloss}{\mathcal{L}_B}
\newcommand{\ls}{\text{LS}^*}

\newcommand{\wa}{\omega_1}
\newcommand{\wb}{\omega_0}

\newcommand{\mlmodel}{\psi}
\newcommand{\trainedmlmodel}{\hat{\psi}}
\newcommand{\fit}{\text{fit}}
\newcommand{\featurizer}{\phi}
\newcommand{\labels}{\mathcal{Y}}

\newcommand{\lsa}{\mathcal{A}_{\text{LS}}}

\newcommand{\spla}{\mathbb{A}_{\text{SPL}}}
\newcommand{\splx}{\mathcal{X}}

\newcommand{\varname}[1]{\texttt{\detokenize{#1}}}

\newcommand{\refeq}[1]{(\ref{#1})}

\newcommand{\lgbma}{LGBMW1}
\newcommand{\lgbmb}{LGBMW-1}

\newcommand{\lns}{LNS}

\newcommand{\setbin}{\mathcal{B}}
\newcommand{\setint}{\mathcal{Q}}
\newcommand{\setcontinuous}{\mathcal{W}}
\newcommand{\lb}{l}
\newcommand{\ub}{u}

\newcommand{\xbbin}{\xb_{\setbin}}
\newcommand{\nsm}{\mathcal{N}}
\newcommand{\state}{\mathcal{S}}
\newcommand{\fixingratio}{r_{\nsm}}
\newcommand{\history}{h}

\newcommand{\keepset}{\mathcal{M}} 

\newcommand{\xblin}{\tilde{\xb}}
\newcommand{\xlin}{\tilde{x}}

\newcommand{\egreedy}{$\varepsilon\text{-greedy}$}
\newcommand{\ucb}{$\alpha\text{-UCB}$}

\newcommand{\xbar}{\bar{x}}
\newcommand{\xbbar}{\bar{\xb}}
\newcommand{\mw}{m_w}

\newcommand{\biga}{\mathcal{A}}
\newcommand{\lpsolver}{\biga_{\text{LP}}}
\newcommand{\bigd}{\mathcal{D}}

\section{Introduction}

The term Large Neighbourhood Search (LNS) was coined by~\cite{shaw_using_1998} to describe his local search heuristic for the Vehicle Routing Problem (VRP). The heuristic leverages a relatedness measure to guide the selection of customers for removal (destroy) and reinsertion (repair). The LNS framework received significant attention for its ability to efficiently navigate in vast solution spaces of large-scale instances. Its flexibility allows for the incorporation of problem-specific knowledge in both the destroy and repair operators. 

Historically, progress on new or complex problems have been limited by the work of domain experts.  Significant research efforts are often required to design effective heuristics for a given problem. General-purpose methodologies such as Machine Learning (ML) have gained popularity for their ability to produce good strategies from data instead of human intuition. In this paper, we study the exploitation of ML for Mixed-Integer Programming (MIP) problems. We are motivated to work on MIPs because they constitute the standard way to model combinatorial optimization (CO) problems found in many areas.  

The research questions this work attempts to answer are as follows: How can supervised learning be integrated within a destroy operator for MIPs and to what extent does it improve the quality of solutions found?  We aim to tackle those questions in a challenging setting where (i) the offline training procedure is minimal and (ii) the training/testing datasets are heterogeneous collections of instances. The goal (i) is important because of scaling laws, which implies that high performance can be achieved by leveraging massive computational resources. This raises concerns of efficiency and fairness when compared to alternative heuristics that require limited or no training. Furthermore, a minimal training procedure is highly accessible and can reasonably be used as a baseline for future research.   The choice (ii) to work with a heterogeneous dataset significantly increases the difficulty of the learning task, but it is required to claim with a high degree of confidence that the approach generalizes well.  We believe that these questions have not been studied before within this constrained setting. Related works on ML-enhanced LNS typically center their experimental setup around one or a few different types of CO problems \citep{sonnerat_learning_2022, huang_searching_2023, liu_machine_2023, feijen_learning-enhanced_2024, robinson_la_rocca_combining_2024}. While some methods demonstrate good performance across different problems, it remains unclear if they can generalize to generic MIPs and remain computationally efficient. 

Our contribution comes from an exhaustive experimental study on the integration of ML tools to solve generic MIPs. We implement a modular LNS matheuristic as a test bench to compare different LNS heuristics. In terms of solution quality, we demonstrate that this matheuristic can significantly improve the performance of the state-of-the-art (SOA) solvers Gurobi~\citep{gurobi} and SCIP~\citep{bestuzheva_scip_2021}. We validate our methodology with noisy oracles which reveal a clear correlation between prediction accuracy and the quality of solutions found. We develop techniques to help improve the ML model with loss adjustments and sampling routines. Our experiment uses the sixth version of the Mixed Integer Programming Library MIPLIB 2017~\citep{gleixner_miplib_2021}. It is a well-recognized dataset of MIP instances aggregated from submissions of many industries and academic institutions. It provides a standardized, curated dataset that represents the variety of MIP problems encountered in practice. Given that we experiment with variable selection strategies that operate on binary variables, we focus on the subset of instances for which the proportion of binary variables is above a given percentage.  

The sections in this paper are organized as follows. Section~\ref{sec:lns:background} presents the problem formulation, the LNS framework and the baseline heuristics. It introduces the concepts and the notation used throughout this article. Section~\ref{sec:lns:related_works} lists related works on the hybridization between ML and CO in the LNS framework. Section~\ref{sec:lns:methodology} describes the methodology we use to integrate supervised learning in a generic LNS heuristic. The experimental results in Section~\ref{sec:lns:results} include the instance selection procedure and performance metrics of various LNS heuristic configurations. Finally, Section~\ref{sec:lns:conclusion} concludes the paper and discusses future work.

\section{Background}  \label{sec:lns:background}
This section begins with the mathematical notation for the problem formulation. This is followed by a description of the LNS framework and a set of baseline LNS heuristics.  

\subsection{Problem Formulation} \label{sec:lns:formulation}

A generic CO problem $\P$ can be stated using the following MIP: 
\begin{align}
    (\P) \quad & \nonumber \\
    \min \quad & c^T \xb \label{eq:mip:obj} \\
    \text{s.t.} \quad & \nonumber \\  
    & A\xb \leq b, \label{eq:mip:constraints} \\
    & x_i \in \{0,1\}, \quad \forall i \in \setbin, \label{eq:mip:bin}\\
    & x_i \in \mathbb{Z} \cap [\lb_i, \ub_i] , \quad \forall i \in \setint, \label{eq:mip:int} \\
    & x_i \in \mathbb{R} \cap [\lb_i, \ub_i] , \quad \forall i \in \setcontinuous. \label{eq:mip:continuous}
\end{align}


The objective function~\refeq{eq:mip:obj} minimizes a linear function which is the scalar product between the coefficients $c$ and the variables $\xb$. Constraints~\refeq{eq:mip:constraints} are represented using the compact notation where $A$ is the matrix of coefficients and $b$ is the right-hand side. Constraints~\refeq{eq:mip:bin}-\refeq{eq:mip:continuous} describe the domains and bounds ($\lb$,  $\ub$) of decision variables partitioned into  $\setbin$, $\mathcal{Q}$ and $\mathcal{W}$, the sets of binary, integer and continuous variables, respectively. The bold formatting is used to denote a vector of variables. 

The continuous relaxation of the problem $\P$, denoted by $\Plin$, is obtained by relaxing the integrality constraints on the variables in sets $\setbin$ and $\setint$. Formally, we can define $\Plin$ as a Linear Program (LP): 

\begin{align}
    (\Plin) \quad & \nonumber \\
    \min \quad & c^T \xb \label{eq:lp:obj} \\
    \text{s.t.} \quad & \nonumber \\ 
    & A\xb \leq b, \label{eq:lp:constraints} \\
    & x_i \in \mathbb{R} , \quad \forall i \in \setbin \cup \setint\cup \setcontinuous,  \label{eq:lp:real}  \\
    & x_i \in [0,1], \quad \forall i \in \setbin, \label{eq:lp:bin}  \\
    & x_i \in [\lb_i, \ub_i] , \quad \forall i \in \setint\cup \setcontinuous. \label{eq:lp:bounds}
\end{align}

The objective function~\refeq{eq:lp:obj} and linear constraints~\refeq{eq:lp:constraints} remain unchanged. Constraints~\refeq{eq:lp:real} state that all variables can take any real value, and constraints~\refeq{eq:lp:bin}-\refeq{eq:lp:bounds} enforce their respective bounds. We use the notation $\xblin$ to refer to the optimal solution of the relaxation $\Plin$. 

In this paper, we assume that the binary variables $\xbbin \in \{x_i| i \in \setbin\}$ represent a non-trivial proportion of the problem $\P$. In Section~\ref{sec:lns:selection}, we reveal the number of instances that satisfy this condition.  

\subsection{LNS Framework}
This section contains our assessment of the LNS framework and sets the related notation for the remaining of this paper. The version we present here is a matheuristic~\citep{fischetti_matheuristics_2018}, but it can also be integrated in the branch-and-bound (\bnb) \citep{hendel_adaptive_2022}. We operate outside the \bnb\ to keep the implementation solver-agnostic and, therefore, more versatile. Algorithm~\ref{alg:lns} describes the LNS procedure. It takes as input a problem $\P$, an initial feasible solution $\xb^0$, a neighbourhood size manager $\nsm$, and an LNS heuristic $\pi$ (also referred to as destroy operator or a variable selection policy). The initial solution is mandatory and can be obtained via a construction heuristic such as the feasibility pump \citep{fischetti_feasibility_2005, bertacco_feasibility_2007} if a feasible solution is difficult to find. After the initialization of the state $\state$, a while loop is defined with $\varname{not_terminated}$, which checks if any of the termination criteria is met. The state keeps track of the various elements of the solution process (elapsed time, solver termination status, iteration count, etc.). In the while loop, $\nsm$ is interrogated to get a fixing ratio $\fixingratio \in [0, 1]$ that represents the proportion of variables to fix in $\P$. To clarify, since a variable can either be fixed or destroyed, one can deduce the number of variables to destroy ($k_d$) and fix ($k_f$) using the following formulas: 
\begin{align}
    &k_d = |\setbin| \cdot (1 - \fixingratio),  \label{eq:lns:fixing} \\
    &k_f = |\setbin| \cdot  \fixingratio. 
\end{align}

The role of the manager $\nsm$ is to add an adaptive neighbourhood size mechanism that dynamically adjusts the sub-problem size. The next step is to build the constraint $\gamma$ using the policy $\pi$, which is then added to $\P$ to create the sub-problem $\P'$ (destroy). In Section~\ref{sec:lns_heuristics}, we cover LNS heuristics that represent specific instantiations of $\pi$.  The sub-problem $\P'$ is then solved with a MIP solver with the incumbent $\xb'$ as a hint (repair). Historical data from the solution process is assigned to $h$ and added to the state $\state$.  As a last step, an $\varname{update_solution}$ routine such as Simulated Annealing~\citep{kirkpatrick_optimization_1983} is called to decide how to update the incumbent. The parameters of the different components used in our experimental setup are explained in Section~\ref{sec:lns:results}. 

\begin{algorithm}[H]
    \caption{Large Neighbourhood Search Matheuristic}\label{alg:lns}
    \begin{algorithmic}[1]
    \Procedure{LNS}{$\P, \xb^0, \nsm, \pi$}

    \State $\xb' \gets \xb^0; \state \gets \state^0$
    \While{$\varname{not_terminated}(\state)$}
        \State $\fixingratio \gets \nsm(\state)$ \Comment{Generate a fixing ratio }
        \State $\gamma \gets \pi(\xb', \fixingratio, \state)$ \Comment{Create a constraint}
        \State $ \P' \gets \varname{update_problem}(\P, \gamma)$ \Comment{Create the sub-problem (destroy)}

        \State $\xb'', \history \gets \varname{solve}(\P', \xb')$ \Comment{Solve the sub-problem (repair)}
        \State $\state \gets \varname{update_state}(\state, \history )$

        \State $\xb' \gets \varname{update_solution}(\xb', \xb'')$ 
    
    \EndWhile
    \State \Return $x'$
    \EndProcedure
    \end{algorithmic} 
\end{algorithm} 

\textbf{Adaptive fixing rate strategy. } In related works (see Section~\ref{sec:lns:related_works}), the fixing ratio is often adjusted dynamically during the solution process. This strategy is referred to as $\nsm$ in Algorithm~\ref{alg:lns}. The proportion of variables fixed in the original problem is critical because the performance of LNS heuristics heavily depends on building neighbourhoods with an appropriate size. One the one hand, fixing more variables makes the sub-problem easier to solve, but it increases the risk of missing good solutions. On the other hand, a low fixing ratio allows for more exploration at increased computational cost. In some cases, we can infer if the neighbourhood size should be increased or decreased based on the termination status of the previous iteration. In particular, the fixing ratio should be decreased if the sub-problem is solved to optimality or is infeasible. We should increase it if the sub-problem is too difficult and no new incumbent is found. In \cite{chmiela_online_2023}, the updated fixing ratio $r^{t+1}$ is 
\begin{align}
    r^{t+1} = (1 \pm \lambda_r ) \cdot r^{t}, \label{eq:updateratio}
\end{align}
where $\lambda_r$ is a scaling parameter and $r^{t}$ is the ratio of the previous iteration. They also impose a lower bound $r_{\text{min}}$ and an upper bound $r_{\text{max}}$ on the value of the ratio.  Other related works use a similar strategy but with a slightly different formula~\citep{hendel_adaptive_2022, huang_local_2023}. 

\subsection{LNS Heuristics}\label{sec:lns_heuristics}

An LNS heuristic typically enforces the following constraint: 
\begin{align}
    \sum_{i \in \keepset} x_i (1 - x'_i) + (1 - x_i) x'_i  = 0, \label{eq:lns:variable_fixing}
\end{align}
which forces $x_i$ to be equal to $x'_i$ for all indices in the fixing set $\keepset \subset \setbin$. The question that an LNS heuristic attempts to answer is: what is a good choice for the set $\keepset$?  One simple answer is to randomly select the variables. The Mutation (or random) heuristic \citep{rothberg_evolutionary_2007} iteratively samples an index without replacement to produce the following fixing set: 
\begin{align}
    \keepset_{Random} = \{i_1, ..., i_{k_f}\}. 
\end{align}
where $i_t$ is the index selected at iteration $t$.

To improve performance, we should fix variables that are likely to have their optimal values in the current solution. To do this, \cite{danna_exploring_2005} proposed Relaxation Induced Neighbourhood Search (RINS). The idea of RINS is to deduce $\keepset$ based on the matching set between the relaxation and the incumbent: 
\begin{align}
    \keepset_{RINS} = \{ i \in \setbin \mid \tilde{x}_i = x_i' \}. 
\end{align}

The Crossover heuristic \citep{rothberg_evolutionary_2007} creates a neighbourhood using the matching set of different feasible solutions. Given solutions $\xb' $ and $\xb''$, this heuristic would fix variables as follows: 
\begin{align}
    \keepset_{Cross} = \{ i \in \setbin \mid  x_i'' = x_i' \}. 
\end{align}

The Local Branching (LB) heuristic differs from previous approaches because the cut generated does not directly fix variables. Instead, LB uses the following soft-fixing constraint: 
\begin{align}\label{eq:lns:lb}
    \sum_{i\in \setbin} x_i  (1 - x_i') + (1 - x_i) x_i' \leq k_{LB}, 
\end{align} 
which restricts the search space within a distance $k_{LB}$ of the incumbent $\xb'$. As discussed in \cite{liu_learning_2022}, the upper bound $k'_{LB}$ on $k_{LB}$ can be calculated using the $L_1$ distance between the incumbent and the LP solution:
\begin{align}\label{eq:lns:kmax}
    k'_{LB} = \| \xb' - \xblin \|_1. 
\end{align} 
An LB cut with a distance parameter above this value will not have any impact on the LP solution. Hence, the right-hand side $k_{LB}$ in~\refeq{eq:lns:lb} is calculated as follows: 
\begin{align}\label{eq:lns:klb}
    k_{LB} = (1 - \fixingratio) k'_{LB}.
\end{align}


Recently, \cite{huang_local_2023} proposed to combine RINS and LB into one heuristic called LB-RELAX. The idea is to add the LB cut~\refeq{eq:lns:lb} to $\Plin$ and produce a fractional solution $\xblin^{LB}$. Then, it computes $\Delta_i = |\xlin^{LB}_i - x_i'|$ the distance between the solution of the relaxation and the incumbent. It identifies the index set of variables for which the distance is zero to build the neighbourhood
\begin{align}
    \keepset_{LB-RELAX} = \{ i \in \setbin \mid  \Delta_i = 0 \}. 
\end{align}
An important consideration we omitted so far is the correction required to ensure that the number of fixed variables is always $k_f$ for all heuristics. If $|\keepset| < k_f$, we need to randomly or greedily \citep{huang_local_2023}  add  variables to $\keepset$. Inversely, if $|\keepset| > k_f$, variables are randomly removed from the set until the target fixing ratio is reached. Note that this correction does not apply to LB and random.

\section{Related Works} \label{sec:lns:related_works}

There is a growing literature on ML-augmented \bnb\ for MIPs~\citep{scavuzzo_machine_2024}. ML predictions can be used to assist various decisions in the \bnb\ framework such as node selection~\citep{khalil_mip-gnn_2022} or branching rules~\citep{khalil_learning_2016}. For our purposes, this section specifically covers related works on the integration of learning within the LNS framework. We start with a description of online methods which is how ML is used in the open-source solver SCIP.  We follow up with recent approaches where learning happens offline. For a broad overview of ML for CO, we refer to the surveys of~\cite{bengio_machine_2021} and \cite{kotary_learning_2021}. 

\textbf{Online methods.} In SCIP 9~\citep{bolusani_scip_2024}, the latest version of the non-commercial solver, self-improving procedures are used to dynamically choose which heuristic to run during the \bnb. Learning in this context occurs online, which means that decisions are taken based on the performance measured on the instance currently being solved. This idea was introduced by~\cite{ropke_adaptive_2006} in the Adaptive Large Neighbourhood Search (ALNS) framework originally applied to the pickup and delivery problem with time windows (PDPTW). The method uses roulette wheel selection based on weights that are periodically adjusted. The goal is to bias the selection towards heuristics that work well on the particular problem. The implementation in SCIP uses modified multi-armed bandit algorithms~\citep{hendel_adaptive_2019, hendel_adaptive_2022, chmiela_online_2023} that compute a reward which considers various factors such as the improvement of the objective value and whether a new solution is found. In~\cite{hendel_adaptive_2022}, they compare \ucb\ (Upper Confidence Bound) and \egreedy\ as different strategies for balancing exploration and exploitation. The strategy \ucb\ explores deterministically and focuses on actions with high uncertainty or potential. They state that \ucb\ achieves the highest solution rate and it is less sensitive to its parameter choice than \egreedy, making it more reliable.

In \cite{robinson_la_rocca_one-shot_2024}, a one-shot learning approach called Probe and Freeze (PNF) is used to guide the search. As a case study, they use the Locomotive Assignment Problem (LAP), which models the flow of locomotives in a railway network. In this problem, the core decisions relate to the assignments of trains to sets of locomotives represented as special ordered sets of type 1 (SOS1). The method tracks intermediate solutions found at the top of the \bnb\ tree to make predictions. This is referred to as the \textit{probing} data. The PNF heuristic fixes a percentage of the SOS1 constraints giving priority to the ones with the lowest entropy. The idea is to quantify uncertainty to limit the risk of compromising the quality of the final solution.  We refer to the work of \cite{ortiz-astorquiza_locomotive_2021} for more details on the problem formulation. 

\textbf{Offline methods.} Most of the recent research efforts in ML have been oriented towards offline learning because it can leverage knowledge stored in large datasets. At inference, the pretrained model is typically called to help guide the search based on what has been learned offline. 

\cite{sonnerat_learning_2022} uses Neural Diving (ND)~\citep{nair_solving_2021} to generate an initial assignment and Neural Neighbourhood Selection (NNS) to guide the search process in LNS. The method uses a bipartite graph representation of MIPs, where nodes represent variables and constraints. The policy uses a Graph Convolutional Network (GCN) to learn good representations.  In ND, the model is trained using an energy function such that assignments with better objective values have higher probabilities. At inference, the predicted distribution is used to generate different sub-problems solved in parallel using MIP solvers. NSS places the variable selection task within a Markov Decision Process (MDP). The GCN is called at every LNS iteration to determine which variables to destroy.  The policy is trained on the optimal assignment obtained by the expert policy LB. On the production packing dataset, the method achieves a final average primal gap that is almost two orders of magnitude smaller than SCIP. 

\cite{liu_learning_2022} use ML to optimize the LB heuristic  in a  two-phase strategy: scaled regression for initial neighbourhood size prediction, followed by reinforcement learning for dynamic adaptation. They also employ Graph Neural Networks (GNN) for the ML model. Unlike previous approaches, it predicts a neighbourhood size at the instance level rather than making variable-level predictions. The best-performing hybrid algorithm (lb-srmrl) significantly outperforms the original LB on set covering problems. 

\cite{huang_searching_2023} propose to use a Contrastive Loss (CL) to train the variable selection policy. This particular loss aims to learn representations where similar samples are close together in the embedding space, while dissimilar samples are far apart. To do so, it leverages both positive and negative examples to provide a richer learning signal. Positive samples come from intermediate solutions that have a relative improvement within a given threshold. For negative samples, they use a controlled degradation technique where perturbations are added incrementally to the optimal set to create nearby solutions that are worse in terms of quality. Their approach, called CL-LNS, uses a Graph Attention Network (GAT) to process the bipartite graph representation. CL-LNS shows the most dramatic improvement on combinatorial auction instances, especially on large instances. For example, at the 60-minute mark, CL-LNS achieves a primal gap of 0.09\%, while the next best method (LB-RELAX) achieves 1.61\%. 

These works are inclined to use GNNs because they are well suited for the data structure of MIPs. For instance, they are invariant to node and edge permutations, they can handle variable-sized inputs, and they can incorporate various node and edge features~\citep{peng_graph_2021}. However, they require the full graph as input and that can create memory capacity concerns for large-scale instances or for training (where batches of graphs are used in the forward pass). While distributed training is a possible solution, it adds implementation complexity and computational costs. To improve scalability, we propose to make predictions on each variable separately. This allows us to run inference iteratively and have a negligible impact on memory usage. \cite{sonnerat_learning_2022} report that training their GNN can take up to 48 hours. In our case, training time is less than a minute, which is a notable distinguishing factor. 

Our methodology takes inspiration from related works that use sampling data as input features for the ML model. \cite{sun_generalization_2021} propose a ML-enhanced reduction strategy for the travelling salesman problem where features are constructed from stochastic sampling of feasible tours. \cite{robinson_la_rocca_combining_2024} adapt this approach to the multicommodity capacitated fixed-charge network design problem and employ slope scaling~\citep{crainic_slope_2004} to generate features. These articles design their approach around heuristics specifically designed for their target problem. Our work aims to generalize every component of the methodology to apply it to generic MIPs.  We provide more details on how we propose to do so in the following section. 


\section{Methodology} \label{sec:lns:methodology}
The methodology section is divided into two parts. The first part explains how we integrate the predictions to create an effective LNS. In the second part, we give instructions related to the training process. 

\subsection{LNS Integration}

We aim to explore how data can help us find a good fixing set $\keepset$ for a given instance. In theory, we can infer one optimal set if we have access to an optimal solution $\xb^*$.  The corresponding fixing set would be 
\begin{align}
    \keepset_{Oracle} = \{ i \in \setbin \mid  x_i' = x_i^* \},  \label{eq:oracle_matching_set}
\end{align}
where variables are fixed to their corresponding value in $\xb^*$. We know that the neighbourhood created by $\keepset_{Oracle}$ will always contain at least one optimal solution. We refer to this method as Oracle LNS (OLNS). In practice, however, it is often intractable to find an optimal solution, so $\xb^*$ refers more broadly to a high-quality solution found by an expert. In the remaining of the paper, it is called the \textit{label solution}. 

We choose to model this problem as a binary classification task. The model $\mlmodel$ is trained in a supervised way to learn to predict the label value $x_i^*$ for each variable. We refer to our method as Supervised LNS (SLNS). At inference, we call $\mlmodel$ on each variable to get an estimation $\xbar_i$ and we use it in a similar fashion to build the corresponding fixing set:
\begin{align}
    \keepset_{\text{SLNS}} = \{ i \in \setbin \mid  x_i' = \xbar_i \}. 
\end{align}

\textbf{Stochastic selection.} The version of SLNS we just introduced is deterministic, which is problematic for two reasons. First, the prediction $\xbbar$ is a noisy estimation of $\xb^*$, and the error rate can be large in practice. Second, the prediction $\xbbar$ would be the same at each iteration and that would lead to poor exploration. Hence, we propose a stochastic version that uses a weighted selection instead: 
\begin{align}
    P(x_i) = \frac{w_i}{\sum_{j \in \setbin } w_j}, 
\end{align}
where the probability $P(x_i)$ to select $x_i$ is proportional to its assigned weight $w_i$. The weight is calculated based on the matching set between $\xb'$ and  $\xbbar$: 

\begin{align}
    w_i = \begin{cases}
        \mw & \text{if } x_i' = \bar{x}_i \\
        1 & \text{otherwise, }
    \end{cases}
\end{align}
where $\mw$ is the weight multiplier and it allows us to tune the distribution. A high $\mw$ brings us closer to deterministic selection and random sampling corresponds to $\mw = 1$. We showcase in Section~\ref{sec:lns:methodology_validation} why this parameterization is important. 

\textbf{Exploration concerns.} The exploration concern is not specific to SLNS as most LNS heuristics can get stuck during the search. This is particularly problematic for LNS heuristics that build the neighbourhood in a deterministic way such as RINS or LB. Inspired by \cite{huang_local_2023}, we add a mechanism in $\pi$ to automatically revert to random selection if no new solution is found after a predetermined number of iterations. 

\subsection{Training Procedure} 
We use supervised learning to train a binary classifier to predict the label value of each variable. We minimize the following binary cross-entropy loss:
\begin{align} \label{eq:lns:binloss}
    \binloss (\xbar, x) = -\wa x \log(\xbar) - \wb (1-x) \log(1-\xbar),
\end{align}
where $x$ is the true label and $\xbar$ is the predicted label. The weights $\wa$ and $\wb$ represent the positive and negative classes, respectively. The weights can be adjusted to control the contribution of each class and improve performance on unbalanced datasets. We measure the effect of weight adjustments in Section~\ref{sec:lns:results}. 

\textbf{Data collection.} Features are built using intermediate solutions found early in the search~\citep{robinson_la_rocca_one-shot_2024}. We call the first phase \textit{probing} (PRB). The probing phase has two main objectives: to find one feasible solution to build a neighbourhood with the LNS heuristic and to collect informative data for the ML model. Probing is done by running the \bnb\ solver for a fixed amount of time and we use its callback to collect feasible and fractional solutions. This strategy is a general and effective way to generate data. However, for large instances, the solver may stay in the root node for a long time, which limits the amount of data we can collect. To mitigate this issue, we propose a second phase called \textit{sampling} (SPL), where we explore the solution space by iteratively solving the relaxation of the problem. 

Algorithm~\ref{alg:lns:spl} describes the SPL procedure. It takes as input the problem $\P$, an initial feasible solution $\xb^0$, and an LP solver $\lpsolver$. The algorithm generates a random fixing ratio $r$ and creates an LB constraint $\gamma_{LB} $ to limit the solution space around the incumbent. The constraint is added to the problem $\Plin$ and solved using the LP solver. The solution $\xbbar$ is then added to the dataset $\bigd$. The SPL procedure is repeated until the time limit is reached. Compared to probing, sampling allows us to collect more data in a shorter amount of time. It avoids the \bnb\ solver, which includes various routines that can be time consuming. We combine the data collected from both phases to train the ML model. 

\begin{algorithm}[H]
    \caption{Sampling procedure for MIPs}\label{alg:lns:spl}
    \begin{algorithmic}[1]
    \Procedure{SPL}{$\P, \xb^0, \lpsolver$}

    \State $\bigd \gets \{ \}$ \Comment{Initialize the dataset}
    \State $\Plin \gets \P$ \Comment{Relax the problem}
    \While{$\varname{not_terminated}$}
        \State $r\sim [0, 1]$ \Comment{Generate a random fixing ratio }
        \State $\gamma_{LB} \gets \pi_{LB}(\xb^0, r)$ \Comment{Create LB constraint}
        \State $ \Plin' \gets \varname{update_problem}(\Plin , \gamma_{LB})$ \Comment{Add LB constraint}
        \State $\xbbar \gets \varname{solve}(\Plin', \lpsolver)$ \Comment{Solve using LP solver}
        \State $\bigd \gets \bigd \cup \{ \xbbar \}$ \Comment{Add the data to the dataset}
    \EndWhile
    \State \Return $\bigd$
    \EndProcedure
    \end{algorithmic} 
\end{algorithm} 

\textbf{Feature engineering.} The number of solutions found during probing and sampling can be large and it varies from one instance to another. To accommodate the fixed-size input requirement of the ML model, a compression step is needed. We use a histogram to represent the distribution of values for each variable. Algorithm~\ref{alg:lns:histogram} describes the procedure to build the histogram data structure from a set of variable values. The input set $X$ contains the values collected for a given variable. The parameter $k_H$ controls the number of bins in the histogram $H$. The algorithm iterates over the set $X$ and assigns each value to a bin. The frequency of each bin represent different input features for the ML model. 

We also include features related to feasible solutions. They are scarcer and more informative than fractional solutions. For that reason, we do not aggregate them in a histogram. Instead, we use the raw values from the 10 best solutions found during probing. 
\begin{algorithm}
    \caption{Procedure to build the histogram data structure from a set of variable values}\label{alg:lns:histogram}
    \begin{algorithmic}[1]
    \Procedure{BuildHistogram}{$X, k_H$}
        \State $n \gets |X|$ \Comment{Number of elements}
        \State $w \gets \frac{1}{k}$ \Comment{Bin width}
        \State $H \gets \text{zero vector of length } k_H$ \Comment{Initialize histogram}
        \For{$i = 1$ to $n$}
            \State $j \gets \lfloor \frac{x_i}{w} \rfloor + 1$ \Comment{Determine bin index}
            \If{$j > k$} 
                \State $j \gets k$ \Comment{Handle edge case for $x_i = 1$}
            \EndIf
            \State $H_j \gets H_j + \frac{1}{n}$ \Comment{Increment frequency}
        \EndFor
        \State \Return $H$
    \EndProcedure
    \end{algorithmic}
\end{algorithm}

\section{Experimental Results} \label{sec:lns:results}
Experimental results are organized into the following subsections. Section~\ref{sec:lns:selection} describes the instance selection process used to build our test set. We study the effect of noisy predictions in Section~\ref{sec:lns:methodology_validation}. In Section~\ref{sec:lns:sampling_results}, we present results related to the sampling routine. Accuracy metrics for our ML model are revealed in Section~\ref{sec:lns:accuracy_metrics}. Section~\ref{sec:lns:main_results} contains the main set of results for the integration of ML predictions in the LNS framework. We include a segmentation analysis in Section~\ref{sec:lns:segmentation} to study our performance on different subsets of instances. Section~\ref{sec:lns:lap_results} presents the last set of results on large LAP instances.

\textbf{Experimental settings. }
Our software stack is written in Julia and we use the JuMP package~\citep{lubin_jump_2023} to interface with the following MIP solvers: Gurobi Optimizer version 10.1~\citep{gurobi}, SCIP version 8~\citep{bestuzheva_scip_2021}, and HiGHS version 1.7~\citep{huangfu_parallelizing_2018}. In the following sections, we run separate experiments for both Gurobi and SCIP, which are the SOA commercial and open-source MIP solvers, respectively. HiGHS is a competitive open-source alternative to SCIP that we use for its strong Interior-Point Method (IPM) implementation (see Section~\ref{sec:lns:sampling_results}). 

By default, the total runtime budget is 10 minutes, probing is 2 minutes, sampling is 1 minute and iteration time is 20 seconds. The time dedicated to search is equal to the time remaining, i.e., 7 minutes, if both probing and sampling are required by the algorithm. We repeat each run three times using a different seed to improve the statistical significance of the results. We include the performance of the MIP solver alone with the full time budget of 10 minutes. The MIP solver is marginally sensitive to the random seed, so we run it using the default seed only. The labels are generated using the best solution found by Gurobi with a 6-hour computational budget. 

In Algorithm~\ref{alg:lns}, the fixing ratio is updated using~\refeq{eq:updateratio} with an initial value of 0.2. The parameter values for $\lambda_r$, $r_{\text{min}}$, and $r_{\text{max}}$ are 1.5, 0.01, and 0.9, respectively.  The $\varname{update_solution}$ function updates the incumbent in a greedy way such that $\xb'$ is always the best known solution. To ensure good exploration, we revert to random LNS if we fail to find an improving solution after two iterations. 

\textbf{Metrics. }
The two primary metrics used to evaluate the performance of algorithms are the primal gap (PG) and the primal integral (PI). The PG is defined as the relative difference between the objective function values of the algorithm's best solution and the best known solution. The PI is calculated as the area under the curve of the primal gap over time. The PG measures the quality of the solutions found, while the PI measures the speed at which the algorithm finds high-quality solutions. Numerical values are often noisy and outliers are common in our experiments. For a more robust evaluation, we compute the mean, the quantiles and the geometric mean (shifted by one unit) of the metrics. The geometric mean is referred to as geomean in the tables of this paper. 

To evaluate the ML model, we report accuracy-related metrics useful for unbalanced datasets. Specifically, we measure the balanced accuracy, which is the average of recall obtained on both classes. In the remainder of the paper, mentions to the accuracy refer to the balanced accuracy. We also report the false negative rate (FNR) and the false positive rate (FPR). The FNR is the proportion of samples where the model predicts 0 instead of 1.  The FPR is the reverse, i.e., the proportion of samples where the model predicts 1 instead of 0.

\subsection{Instance Selection} \label{sec:lns:selection}
This section describes the process of choosing which instance from MIPLIB to consider.  Based on data-driven analyses, the paper on MIPLIB~\citep{gleixner_miplib_2021} tags a subset of 240 instances as \textit{benchmark suitable}. However, it is designed to benchmark \bnb\ solvers where the main goal is to prove optimality within a 2-hour window. For our purposes, these instances are often too easy. To increase our dataset size, we consider the full collection of 1,065 instances. The following is the list of conditions \mkbluethesis{that} each instance must satisfy:

\begin{enumerate}
    \item The instance must be feasible. \label{item:feasible}
    \item The instance must not generate any error. \label{item:conditioned}
    \item The instance must contain a reasonable proportion of binary variables. \label{item:binary}
    \item A feasible solution must be available.  \label{item:solution}
    \item The instance should not be trivial to solve. \label{item:nontrivial}
\end{enumerate}

The condition~\refeq{item:feasible} can be determined using tags already available. We observe that some instances are tagged as \varname{no_solution} or \varname{infeasible}. Condition~\refeq{item:conditioned} refers to a collection of criteria such as out of memory, loading or numerical errors that make the instance impractical to work with. Condition~\refeq{item:binary} is required because we apply the LNS framework on binary variables and, thus, we need a condition to ensure they are present in our test instances. We filter out instances where binary variables represent less than 10\% of all variables. Table~\ref{tab:instance_selection} contains how many instances satisfy each condition explained so far. We are left we 791 instances after this first filtering step. 

 \begin{table}[h]
\centering
\caption{Instance selection based on the first three conditions}
\label{tab:instance_selection}
\begin{tabular}{lcc}
\toprule
{Selection} & {Count} \\
\midrule
Total & 1065 \\
Error & 24 \\
Infeasible & 60 \\
Low binary ratio & 190 \\
Selected & 791 \\
\bottomrule
\end{tabular}
\end{table}

 The second selection phase relates to conditions~\refeq{item:solution} and \refeq{item:nontrivial}. Unlike before, these require intermediary experiments. Condition~\refeq{item:solution} is a requirement of the LNS framework itself. To find a solution, we run each generic solver for 2 minutes with default settings. Table~\ref{tab:sol_node_count_distribution} includes \mkbluethesis{statistics for metrics related to the solution and node counts. For both metrics, we report the corresponding average and the number of instances for which the metric is exactly equal to 0 or 1.} To clarify, the node count corresponds to the number of fractional solutions that we were able to collect via the callback of the solver. Results reveal that Gurobi generates significantly more solutions compared to SCIP (around 2x for feasible solutions, and 75x for fractional solutions). The insight relevant to the selection process is that Gurobi and SCIP fail to find a feasible solution for 106 and 220 instances, respectively. 

\begin{table}[h]
\centering
\caption{Solution and node count distribution for 2-minute runs}
\label{tab:sol_node_count_distribution}
\begin{tabular}{lccccccc}
\toprule
{Solver} & \multicolumn{3}{c}{Feasible solutions} & \multicolumn{3}{c}{B\&B nodes} \\
{} & {0} & {1} & {mean} & {0} & {1} & {mean} \\
\midrule
GRB & 106 & 76 & 8.06 & 138 & 37 & 3316.73 \\
SCIP & 220 & 127 & 4.50 & 227 & 125 & 44.42 \\
\bottomrule
\end{tabular}
\end{table}

 Condition~\refeq{item:nontrivial} states that the instance must be sufficiently difficult to justify the use of a primal heuristic. If an exact solver can close the optimality gap relatively quickly, the addition of LNS heuristics is unlikely to improve performance~\citep{chmiela_online_2023}.  We use the optimality gap found after a 10-minute run to measure the difficulty of an instance and we report that metric in Table~\ref{tab:optimality_gap_distribution}.   We note that 171 and 258 instances have an optimality gap above 10\% after 10 minutes for Gurobi and SCIP, respectively. Gurobi reaches near-optimal solutions for most instances.

\begin{table}[h]
\centering
\caption{Optimality gap distribution for 10-minute runs}
\label{tab:optimality_gap_distribution}
\begin{tabular}{lccc}
\toprule
{Optimality gap} & {GRB} & {SCIP} \\
\midrule
Below 0.1\% & 442 & 273 \\
Between 0.1\% and 10\% & 92 & 112 \\
Above 10\% & 171 & 258 \\
No solution & 86 & 148 \\
Total & 791 & 791 \\
\bottomrule
\end{tabular}
\end{table}

Table~\ref{tab:instance_selection_from_results} summarizes the selection process with the number of instances related to conditions~\refeq{item:solution} and \refeq{item:nontrivial} for both solvers. An instance is selected if it satisfies these two conditions for each solver. We end up with 126 instances for the test set that will be used to benchmark different LNS heuristics in the following sections. The tag distribution for the test set is provided in Table~\ref{tab:tags_for_selected}. \mkbluethesis{Each MIPLIB instance has a set of publicly available tags that typically indicate the problem type.}

\begin{table}[h]
\centering
\caption{Instance selection based on results from 2-minute and 10-minute runs.}
\label{tab:instance_selection_from_results}
\begin{tabular}{lcc}
\toprule
{Selection} & {Count} \\
\midrule
Solution found with GRB with 2-minute run (condition (4)) & 685 \\
Solution found with SCIP with 2-minute run (condition (4)) & 571 \\
Hard for GRB with 10-minute run (condition (5)) & 171 \\
Hard for SCIP with 10-minute run (condition (5)) & 258 \\
Intersection between conditions (4) and (5) for both solvers & 126 \\
\bottomrule
\end{tabular}
\end{table}

\begin{table}[h]
\centering
\caption{Tags distribution for test set of 126 instances}
\label{tab:tags_for_selected}
\begin{tabular}{lcc}
\toprule
{Tag} & {Count} \\
\midrule
mixed\_binary & 71 \\
variable\_bound & 67 \\
precedence & 64 \\
invariant\_knapsack & 49 \\
set\_partitioning & 48 \\
binary & 39 \\
general\_linear & 35 \\
set\_packing & 34 \\
benchmark\_suitable & 33 \\
aggregations & 31 \\
cardinality & 27 \\
decomposition & 25 \\
knapsack & 20 \\
set\_covering & 20 \\
benchmark & 19 \\
numerics & 18 \\
binpacking & 8 \\
integer\_knapsack & 4 \\
equation\_knapsack & 3 \\
indicator & 2 \\
\bottomrule
\end{tabular}
\end{table}

\subsection{Methodology Validation}\label{sec:lns:methodology_validation} 
This section aims to demonstrate that our methodology is sound. To do so, we experiment with OLNS, an algorithm that uses the label solution to guide the search. We believe that this experiment is crucial to validate our experimental setup. It sets a lower bound on the primal gap that we can expect with our approach. Also, it allows us to simulate the performance we could obtain with supervised learning. By adding noise to the oracle, we imitate the behaviour of the noisy predictions given by a ML model. We also include a deterministic version (DOLNS) that selects variables based on the matching set~\refeq{eq:oracle_matching_set}. We run each configuration with three different seeds on 10 randomly selected instances from the test set. 

\begin{table}[h]
\centering
\caption{Shifted geometric mean of the primal gap for different oracle configurations}
\label{tab:oracle_error_rates}
\begin{tabular}{lccccc}
\toprule
{Oracles} & \multicolumn{4}{c}{Error rate} \\
{} & {0.0} & {0.1} & {0.3} & {0.5} \\
\midrule
OLNS-2 & 3.06 & 3.71 & 3.65 & 3.93 \\
OLNS-100 & 1.43 & 2.73 & 3.74 & 4.05 \\
OLNS-1000 & 1.39 & 2.65 & 3.69 & 3.99 \\
OLNS-10000 & 1.33 & 2.71 & 3.93 & 4.16 \\
DOLNS & 1.07 & 3.89 & 3.82 & 4.03 \\
\bottomrule
\end{tabular}
\end{table}

Table~\ref{tab:oracle_error_rates} reveals significant trends under varying noise conditions. In the absence of noise (0.0 error rate), DOLNS demonstrates the best performance with the lowest primal gap (1.07). However, there is a steep performance degradation when any noise is added. Unlike DOLNS, OLNS is more resilient to noise. OLNS-2 maintains the most consistent performance across error rates. Algorithms with higher weight multipliers (OLNS-100, OLNS-1000, OLNS-10000) perform well at lower error rates but show a notable decrease in performance as noise levels rise. These findings suggest that the selection of a weight multiplier should be based on the anticipated error rate of the ML predictions. Given the difficulty of our task, we expect a relatively high prediction error. Hence, we use a weight multiplier of 2 in the following experiments. 

These results are meaningful because they validate our methodology. They demonstrate a clear correlation between our target performance (primal gap) and the prediction accuracy. This implies that the solution process for generic MIPs can be improved by increasing the accuracy of a supervised learning model. Since advances in ML happen quickly, we expect that our ability to predict high-quality solutions will improve over time. As accuracy increases, we must increase the weight multiplier to take advantage of it. 


\subsection{Sampling Results} \label{sec:lns:sampling_results}
The goal of the sampling routine is to generate additional data to help improve the ML model predictions. To do so, we iteratively solve the continuous relaxation with different LB cuts as described in Algorithm~\ref{alg:lns:spl}. We give a time limit of 5 seconds to the LP solver, but it is not always sufficient to find the optimal solution within the default tolerance.   The simplex algorithm~\citep{dantzig1951maximization} is typically the default LP solver in \bnb\ solvers because it can easily be warm-started. However, it suffers from exponential worst-case time complexity and it can struggle to converge for some instances. We noticed in Table~\ref{tab:sol_node_count_distribution} \mkbluethesis{the} poor performance with SCIP which uses an LP solver primarily designed for simplex-based methods (SoPlex). To improve the performance of the sampling routine, we experiment with IPM \citep{karmarkar_new_1984}. This method offers a polynomial-time complexity and it can be more effective for large-scale problems. We use the high-quality implementation of IPM in HiGHS~\citep{huangfu_parallelizing_2018}. We run sampling with three different LP solver configurations available in HiGHS: simplex, IPM, and Primal-Dual LP (PDLP). We adjust the following gap tolerances to speed up the LP solvers: $\varname{pdlp_d_gap_tol} = 10^{-3}$, 
$\varname{ipm_optimality_tolerance} = 10^{-4}$. 

\begin{table}[h]
\centering
\caption{Instance count at different sample sizes for SPL runs with 1-minute time limit}
\label{tab:sample_size_distribution}
\begin{tabular}{lccccc}
\toprule
{Sample size} & {GRB} & \multicolumn{3}{c}{HiGHS} \\
{} & {default} & {ipm} & {simplex} & {pdlp} \\
\midrule
$[0, 10]$ & 39 & 23 & 12 & 63 \\
$[11, 20]$ & 9 & 47 & 44 & 25 \\
$[21, 30]$ & 7 & 6 & 8 & 4 \\
$[31, 50]$ & 71 & 50 & 62 & 34 \\
\bottomrule
\end{tabular}
\end{table}

Table~\ref{tab:sample_size_distribution} contains the instance count for various sample sizes generated using Gurobi and HiGHS. For all LP solvers, instances are unevenly distributed across the sample size brackets. A good proportion of instances have fewer than 10 samples generated because the time allocated to sample is relatively low. For the largest instances, we suspect that LP solvers do not have enough time to converge within the 5-second time limit. However, this limit is necessary to allow for appropriate exploration during the 1-minute time budget. 

\subsection{ML Model and Accuracy Metrics} \label{sec:lns:accuracy_metrics}
We discussed in Section~\ref{sec:lns:methodology_validation} how the prediction accuracy relates to the quality of the solutions found via a noisy oracle. In this section, we discuss the ML model choice and evaluate its accuracy-related metrics on the test set. 

\textbf{Motivations for our ML model.} We use LightGBM (LGBM) ~\citep{ke_lightgbm_2017}, a boosted tree library, to fit our data and produce a binary classifier. Surveys that compare the performance of various ML models demonstrate that boosted trees are competitive on a wide range of tasks~\citep{caruana_empirical_2006, borisov_deep_2022_2}. They are also more straightforward and cheaper to use than neural networks. Thus, we believe LGBM is a reasonable baseline given our goal of keeping the training procedure minimal. 

\textbf{Dataset splits.}  The established data splitting strategy is to keep around 20\% of samples for testing and the remaining for training. However, given the limited number of instances available, we believe that this approach is not appropriate. Instead, we propose to adjust the training set for each instance to allow us to test on all of them (126 in total). Specifically, the training set always contains the 125 remaining instances after removing the test instance itself. We can afford to train a different model for each instance because LGBM requires less than a minute to train. 

\textbf{Bias adjustments. } In Section~\ref{sec:lns:methodology}, we introduced the weighted version of the binary loss~\refeq{eq:lns:binloss}. This loss includes two weights that control the relative importance of both classes. By default, LGBM uses equal weights and that causes issues when the class distribution is unbalanced. In MIPs found in MIPLIB, the label 0 is more likely than 1 for any given binary variable.  Often, there is a significant asymmetry, where less than 1\% of binary variables have the value 1 in the label solution. In this context, the default LGBM configuration will systematically learn to predict 0. We utilize bias adjustments in the loss to help alleviate this issue. We experiment with three different bias adjustments: W1 ($\wb = 0.25$, $\wa=0.75$), W2 ($\wb = 0.10$, $\wa=0.90$), W3 ($\wb = 0.05$, $\wa=0.95$).

\textbf{Discussion. } Table~\ref{tab:accuracy_metrics} shows the metrics evaluated on the test set: accuracy (Acc), false negative rate (FNR) and false positive rate (FPR). The first suffix (W1 - W3) refers to bias adjustments and the second suffix refers to the data source used for feature engineering. We can highlight three takeaways from this table. First, LGBM performs better with the features collected by Gurobi. Since Gurobi finds better quality solutions during probing, LGBM can infer more accurately what the label solution is. Second, unbiased LGBM has the best accuracy across both solvers, but it also comes with a high FNR. As we increase the bias, LGBM's FNR decreases systematically, which aligns with our expectations. This comes at the cost of a higher FPR. We will discuss the effect of that bias at inference in Section~\ref{sec:lns:main_results}.  The third takeaway relates to the effect of the data used for feature engineering. For SCIP, the relative accuracy improvement from SPL is significant (data from HiGHS). The accuracy jumps from 35\% to 81\% when sampling features are added. This allows us to claim that the sampling routine we propose can extract information that is useful for our prediction task. However, the impact of sampling features is marginal for Gurobi. This disparity between solvers can be explained by the results in Table~\ref{tab:sol_node_count_distribution}. Notice how Gurobi and SCIP explore 3313.14 and 44.34 nodes on average, respectively. The amount of data collected from probing with SCIP is relatively low and that is why it benefits from sampling disproportionally. 

\begin{table}[h]
\centering
\caption{Accuracy related metrics for different bias adjustments}
\label{tab:accuracy_metrics}
\begin{tabular}{lccccccc}
\toprule
{Model} & \multicolumn{3}{c}{GRB} & \multicolumn{3}{c}{SCIP} \\
{} & {Acc} & {FNR} & {FPR} & {Acc} & {FNR} & {FPR} \\
\midrule
LGBMW3-PRB & 0.27 & 0.02 & 0.91 & 0.22 & 0.02 & 0.92 \\
LGBMW2-PRB & 0.36 & 0.06 & 0.79 & 0.22 & 0.02 & 0.91 \\
LGBMW1-PRB & 0.75 & 0.33 & 0.30 & 0.24 & 0.03 & 0.89 \\
LGBM-PRB & 0.84 & 0.59 & 0.13 & 0.35 & 0.17 & 0.74 \\
LGBMW3-SPL & 0.19 & 0.02 & 0.97 & 0.18 & 0.03 & 0.93 \\
LGBMW2-SPL & 0.32 & 0.05 & 0.83 & 0.20 & 0.06 & 0.90 \\
LGBMW1-SPL & 0.77 & 0.34 & 0.27 & 0.51 & 0.36 & 0.53 \\
LGBM-SPL & 0.85 & 0.55 & 0.13 & 0.81 & 0.76 & 0.15 \\
\bottomrule
\end{tabular}
\end{table}

\subsection{Main Results}\label{sec:lns:main_results}

In this section, we compare various LNS configurations on the set of instances selected in Section~\ref{sec:lns:selection}.  Our comparison includes four LNS heuristic baselines that we introduced in Section~\ref{sec:lns_heuristics}: random, LB, LB-Relax, and RINS. In the next paragraphs, we discuss the results using Gurobi and SCIP, separately. Note that in the following discussion, mentions to the average refers to the geometric average which is more appropriate due to the presence of outliers. 

\textbf{Gurobi results.} In Tables~\ref{tab:lns:primal_gap_grb} and \ref{tab:lns:primal_integral_grb}, we show the primal gap and primal integral, respectively. The first three rows of each table represent different Gurobi configurations. Gurobi with $\varname{MIPFocus=1}$ (grb-f1) increases the computational budget given to heuristics. The configuration grb-f1-t3 has the ability to use three threads. Across the three Gurobi configurations, the default one (grb) has the worst primal gap on average but it has the most wins. The primal integral in multi-threaded environment is not available because callbacks require a single-threaded process. 

Within the set of baselines, the random LNS heuristic has the most wins followed by RINS, LB-Relax and LB in that order. Random LNS performs surprisingly well both in terms of primal gap (3.96) and primal integral (8.14). It also has the most wins for the two metrics. Interestingly, we find that RINS outperforms LB-Relax. This is caused by the computational cost of solving the continuous relaxation at every iteration in LB-Relax. The instances in our test set are large and the relaxation is often difficult to solve. RINS only needs the relaxation of the original problem. The LB heuristic has the worst performance overall because it uses a soft-fixing constraint that creates relatively large sub-problems. 

Our heuristic is evaluated using the four LGBM configurations introduced in Section~\ref{sec:lns:accuracy_metrics}. The best version is SLNS-LGBMW2-PRB. It uses the probing data and the second bias adjustment. It has the second-best performance overall, after random, with a primal gap of 4.45 and primal integral of 8.18. The SLNS heuristic with the predictions of unbiased LGBM performs poorly on average. These results confirm that the bias adjustments are useful to improve performance at inference. The sampling data does not produce a meaningful effect, which is consistent with the results in Table~\ref{tab:accuracy_metrics}. 

\textbf{SCIP results.} We present the primal gap and integral in Tables~\ref{tab:lns:primal_gap_scip} and \ref{tab:lns:primal_integral_scip}, respectively. Similar to Gurobi, SCIP has a parameter that sets heuristics emphasis as \textit{aggressive}. This enables improvement heuristics within the \bnb\ of SCIP. In version 9 \citep{bolusani_scip_2024}, a new scheduler for LNS and diving heuristics have been integrated in the solver. For completeness, we also include SCIP 9 results which we run separately because it is not currently supported by the Julia package. All other experiments use version 8. We observe that the latest version upgrade offers a meaningful improvement from 15.98 to 14.64 for the average primal gap. With SCIP, the solver-agnostic callback only returns fractional solutions. That is why the primal integral does not appear for the solver alone. In general, the primal integral is strongly correlated with the primal gap because we use a small total time budget of 10 minutes. 

The relative performance of the four baselines is the same as Gurobi. Random LNS remains the best performing option with an average gap of 6.26 and 24 wins. For SLNS, the sampling data is more useful compared to Gurobi. Its best configuration is SLNS-LGBMW1-SPL with 21 wins and 7.36 for the average primal gap. 

\begin{table}[h]
\centering
\caption{Primal gap (\%) for different LNS heuristics with Gurobi}
\label{tab:lns:primal_gap_grb}
\begin{tabular}{lccccccc}
\toprule
{Scenario} & \multicolumn{3}{c}{Quantiles} & {Mean} & {Geomean} & {Wins} \\
{} & {0.1} & {0.5} & {0.9} & {} & {} & {} \\
\midrule
grb & 0.00 & 2.68 & 88.83 & 134.71 & 6.12 & 20 \\
grb-f1 & 0.00 & 2.51 & 85.22 & 131.82 & 5.73 & 12 \\
grb-f1-t3 & 0.00 & 2.29 & 84.33 & 23.71 & 5.14 & 19 \\
LB & 0.00 & 6.67 & 159.22 & 257.89 & 11.06 & 12 \\
LB-Relax & 0.00 & 2.50 & 57.17 & 44.47 & 5.77 & 21 \\
random & 0.00 & 1.82 & 30.30 & 15.04 & 3.91 & 36 \\
RINS & 0.00 & 2.91 & 44.85 & 18.48 & 4.93 & 25 \\
SLNS-LGBMW3-PRB & 0.00 & 1.98 & 35.59 & 24.60 & 4.65 & 24 \\
SLNS-LGBMW2-PRB & 0.00 & 1.91 & 39.08 & 19.19 & 4.47 & 27 \\
SLNS-LGBMW1-PRB & 0.00 & 1.98 & 50.64 & 22.37 & 4.79 & 26 \\
SLNS-LGBM-PRB & 0.00 & 2.63 & 58.78 & 22.22 & 5.33 & 21 \\
SLNS-LGBMW3-SPL & 0.00 & 1.96 & 40.64 & 21.21 & 4.68 & 28 \\
SLNS-LGBMW2-SPL & 0.00 & 1.93 & 42.32 & 21.80 & 4.72 & 23 \\
SLNS-LGBMW1-SPL & 0.00 & 2.41 & 54.84 & 21.50 & 5.15 & 19 \\
SLNS-LGBM-SPL & 0.00 & 2.49 & 59.81 & 23.71 & 5.39 & 22 \\
\bottomrule
\end{tabular}
\end{table}

\begin{table}[h]
\centering
\caption{Primal integral for different LNS heuristics with Gurobi}
\label{tab:lns:primal_integral_grb}
\begin{tabular}{lccccccc}
\toprule
{Scenario} & \multicolumn{3}{c}{Quantiles} & {Mean} & {Geomean} & {Wins} \\
{} & {0.1} & {0.5} & {0.9} & {} & {} & {} \\
\midrule
grb & 0.43 & 10.30 & 89.80 & 25.47 & 10.74 & 10 \\
grb-f1 & 0.53 & 9.52 & 89.52 & 24.75 & 10.48 & 3 \\
grb-f1-t3 & nan & nan & nan & nan & nan & 0 \\
LB & 0.03 & 7.52 & 94.74 & 28.22 & 9.80 & 12 \\
LB-Relax & 0.07 & 8.07 & 79.44 & 22.59 & 9.45 & 5 \\
random & 0.04 & 6.82 & 54.07 & 18.22 & 8.05 & 31 \\
RINS & 0.04 & 7.69 & 63.79 & 20.29 & 8.79 & 9 \\
SLNS-LGBMW3-PRB & 0.04 & 7.09 & 56.87 & 19.41 & 8.28 & 9 \\
SLNS-LGBMW2-PRB & 0.04 & 6.72 & 65.17 & 19.28 & 8.15 & 16 \\
SLNS-LGBMW1-PRB & 0.06 & 7.45 & 76.00 & 20.98 & 8.61 & 13 \\
SLNS-LGBM-PRB & 0.04 & 7.05 & 72.97 & 21.23 & 8.84 & 10 \\
SLNS-LGBMW3-SPL & 0.03 & 6.83 & 68.72 & 19.76 & 8.20 & 12 \\
SLNS-LGBMW2-SPL & 0.04 & 6.89 & 65.38 & 19.78 & 8.30 & 9 \\
SLNS-LGBMW1-SPL & 0.05 & 6.41 & 72.51 & 20.45 & 8.28 & 10 \\
SLNS-LGBM-SPL & 0.04 & 6.52 & 85.26 & 21.45 & 8.63 & 16 \\
\bottomrule
\end{tabular}
\end{table}

\begin{table}[h]
\centering
\caption{Primal gap (\%) for different LNS heuristics using SCIP}
\label{tab:lns:primal_gap_scip}
\begin{tabular}{lccccccc}
\toprule
{Scenario} & \multicolumn{3}{c}{Quantiles} & {Mean} & {Geomean} & {Wins} \\
{} & {0.1} & {0.5} & {0.9} & {} & {} & {} \\
\midrule
LB & 0.14 & 22.24 & 498.12 & 3895.16 & 26.77 & 3 \\
LB-Relax & 0.00 & 8.84 & 196.15 & 3326.96 & 12.86 & 10 \\
random & 0.00 & 3.26 & 73.59 & 3090.05 & 6.26 & 24 \\
RINS & 0.00 & 6.25 & 75.69 & 3101.75 & 8.82 & 17 \\
SLNS-LGBMW3-PRB & 0.00 & 4.38 & 126.28 & 3122.55 & 7.49 & 15 \\
SLNS-LGBMW2-PRB & 0.00 & 4.46 & 123.69 & 2935.78 & 7.66 & 16 \\
SLNS-LGBMW1-PRB & 0.00 & 4.42 & 108.57 & 3118.66 & 7.63 & 15 \\
SLNS-LGBM-PRB & 0.00 & 4.11 & 109.64 & 2871.21 & 7.63 & 16 \\
SLNS-LGBMW3-SPL & 0.00 & 5.06 & 106.48 & 3105.97 & 7.99 & 12 \\
SLNS-LGBMW2-SPL & 0.00 & 4.88 & 120.22 & 3105.74 & 7.76 & 13 \\
SLNS-LGBMW1-SPL & 0.00 & 4.38 & 100.07 & 3088.26 & 7.36 & 21 \\
SLNS-LGBM-SPL & 0.00 & 5.72 & 148.95 & 3041.16 & 9.12 & 10 \\
scip & 0.00 & 17.51 & 239.26 & 248.18 & 19.10 & 6 \\
scip-9-aggressive & 0.00 & 13.09 & 175.00 & 244.66 & 14.64 & 11 \\
scip-aggressive & 0.00 & 16.04 & 243.81 & 245.54 & 15.98 & 9 \\
\bottomrule
\end{tabular}
\end{table}

\begin{table}[h]
\centering
\caption{Primal integral for different LNS heuristics using SCIP}
\label{tab:lns:primal_integral_scip}
\begin{tabular}{lccccccc}
\toprule
{Scenario} & \multicolumn{3}{c}{Quantiles} & {Mean} & {Geomean} & {Wins} \\
{} & {0.1} & {0.5} & {0.9} & {} & {} & {} \\
\midrule
LB & 0.65 & 25.03 & 96.35 & 40.07 & 19.28 & 3 \\
LB-Relax & 0.42 & 17.19 & 95.44 & 33.44 & 15.27 & 8 \\
random & 0.45 & 10.70 & 89.60 & 24.33 & 11.08 & 30 \\
RINS & 0.34 & 18.25 & 85.42 & 28.63 & 13.64 & 12 \\
SLNS-LGBMW3-PRB & 0.23 & 11.88 & 96.82 & 26.56 & 11.50 & 11 \\
SLNS-LGBMW2-PRB & 0.23 & 12.16 & 97.16 & 26.97 & 11.62 & 10 \\
SLNS-LGBMW1-PRB & 0.26 & 13.06 & 96.64 & 26.45 & 11.68 & 12 \\
SLNS-LGBM-PRB & 0.35 & 12.12 & 96.74 & 26.91 & 11.85 & 7 \\
SLNS-LGBMW3-SPL & 0.27 & 11.72 & 96.61 & 26.14 & 11.54 & 14 \\
SLNS-LGBMW2-SPL & 0.28 & 10.93 & 97.23 & 26.13 & 11.47 & 18 \\
SLNS-LGBMW1-SPL & 0.28 & 13.18 & 96.72 & 25.02 & 11.34 & 17 \\
SLNS-LGBM-SPL & 0.47 & 13.54 & 96.90 & 27.30 & 12.22 & 9 \\
\bottomrule
\end{tabular}
\end{table}

\textbf{Discussion. } 
The most noteworthy takeaway is the relative performance of the random LNS. This has multiple implications. For an arbitrary MIP, the best option is to select variables at random. This is compatible with the generally accepted idea that heuristics must exploit the structure of a problem to be effective. Instances in MIPLIB often have unique structures and they are loosely related to each other. Thus, there is not an obvious structure to exploit in this dataset. 

The effectiveness of the random LNS can be explained by its unbiased nature. At every iteration, a completely new and unrelated neighbourhood is created, which helps exploration.  All other heuristics are biased in some way, and that can lead to systematically exploring regions of the solution space that contain sub-optimal solutions. RINS and LB-Relax, for example, bias the search towards solutions near the continuous relaxation of the problem. This strategy will yield poor performance if near-optimal solutions are far from the relaxation.  Our approach, SLNS, bias the search based on the predictions obtained by the ML model. Our predictions are often inaccurate, which leads us to search in the wrong regions of the solution space. 

The silver lining is that our experiments on SLNS reveal that we are in the right direction. At the forefront, we outperform all other non-random heuristics. Also, it is clear that weight adjustments in the binary loss help improve performance. It is generally a good idea to adjust the loss to reduce the FNR. Next, the addition of sampling features can be beneficial in cases where the probing data is lacking (for SCIP in particular). Another notable outcome of our research is the demonstration that the LNS matheuristic (Algorithm~\ref{alg:lns}) can significantly improve the discovery of high-quality solutions in MIPs compared to SOA \bnb\ solvers.  The performance benefit does not come from a sophisticated variable selection strategy, but rather, it is the result of breaking down the problem in sub-problems that are easier to handle.


\subsection{Segmentation Analysis} \label{sec:lns:segmentation} 
In the previous section, we presented results aggregated over all instances and found that the random LNS performs best on average. However, this analysis ignores the underlying distribution of instances. As discussed in Section~\ref{sec:lns:selection}, MIPLIB contains clusters of instances that share common characteristics (see Table~\ref{tab:tags_for_selected}). In the following, we discuss an analysis on various segments of the dataset. 

Tables~\ref{tab:segmentation_grb} and \ref{tab:segmentation_scip} summarize the segmentation analysis for Gurobi and SCIP, respectively. Each row corresponds to a different segment of the population of instances. The first four rows refer to segments relative to the corresponding median.  For example, the high binary segment represents the set of instances where the number of binary variables is higher than the median value across all instances. The remaining segments refer to instances with the corresponding MIPLIB tag. Results with Gurobi reveal that SLNS is better than random on two segments: high sample size and low binary. Note that these two categories are related because it is easier to produce samples on smaller instances. On the high sample size segment, the average geomean is 2.79 with SLNS and 3.01 for random.

This pattern does not appear with SCIP, because features are primarily generated using the data from sampling. Unlike the probing data, the sampling data is highly correlated with the incumbent, and that can explain why more data does not improve performance in this case. There is a mismatch between the mean and the geomean on SCIP's data because it contains many outliers. In most tag categories, results suggest that SLNS is more resilient to outliers compared to random. 

\begin{table}[h]
\centering
\caption{Segmentation analysis with Gurobi (bold underscore for the best values)}
\label{tab:segmentation_grb}
\begin{tabular}{lcccccc}
\toprule
{Segment} & {Count} & \multicolumn{2}{c}{random} & \multicolumn{2}{c}{SLNS-LGBMW2-PRB} \\
{} & {} & {Mean} & {Geomean} & {Mean} & {Geomean} \\
\midrule
High binary & 63 & \underline{\textbf{23.43}} & \underline{\textbf{5.87}} & 32.51 & 7.75 \\
Low binary & 63 & 7.26 & 2.74 & \underline{\textbf{6.53}} & \underline{\textbf{2.71}} \\
High sample size & 63 & 7.73 & 3.01 & \underline{\textbf{7.25}} & \underline{\textbf{2.79}} \\
Low sample size & 63 & \underline{\textbf{22.96}} & \underline{\textbf{5.34}} & 31.66 & 7.49 \\
Binpacking & 8 & 52.04 & \underline{\textbf{7.88}} & \underline{\textbf{48.95}} & 8.34 \\
Set packing & 34 & \underline{\textbf{26.58}} & \underline{\textbf{3.76}} & 34.48 & 5.09 \\
Invariant knapsack & 49 & \underline{\textbf{23.25}} & \underline{\textbf{5.39}} & 29.09 & 6.36 \\
Knapsack & 59 & \underline{\textbf{21.44}} & \underline{\textbf{5.11}} & 26.61 & 5.83 \\
Mixed binary & 71 & \underline{\textbf{17.27}} & \underline{\textbf{4.58}} & 23.86 & 5.50 \\
Benchmark suitable & 33 & \underline{\textbf{14.2}} & \underline{\textbf{2.09}} & 20.91 & 2.35 \\
\bottomrule
\end{tabular}
\end{table}

\begin{table}[h]
\centering
\caption{Segmentation analysis with SCIP (bold underscore for the best values)}
\label{tab:segmentation_scip}
\begin{tabular}{lcccccc}
\toprule
{Segment} & {Count} & \multicolumn{2}{c}{random} & \multicolumn{2}{c}{SLNS-LGBMW1-SPL} \\
{} & {} & {Mean} & {Geomean} & {Mean} & {Geomean} \\
\midrule
High binary & 63 & 6161.57 & \underline{\textbf{10.35}} & \underline{\textbf{6158.76}} & 13.94 \\
Low binary & 63 & 18.53 & \underline{\textbf{3.79}} & \underline{\textbf{17.76}} & 3.89 \\
High sample size & 63 & \underline{\textbf{9.79}} & \underline{\textbf{3.35}} & 9.97 & 3.64 \\
Low sample size & 63 & 6170.32 & \underline{\textbf{11.7}} & \underline{\textbf{6166.55}} & 14.88 \\
Binpacking & 8 & \underline{\textbf{28.31}} & \underline{\textbf{10.65}} & 46.56 & 12.06 \\
Set packing & 34 & 105.49 & \underline{\textbf{9.11}} & \underline{\textbf{88.78}} & 10.88 \\
Invariant knapsack & 49 & 77.27 & \underline{\textbf{7.62}} & \underline{\textbf{69.4}} & 10.06 \\
Knapsack & 59 & 68.76 & \underline{\textbf{6.96}} & \underline{\textbf{61.11}} & 8.87 \\
Mixed binary & 71 & 59.51 & \underline{\textbf{6.35}} & \underline{\textbf{53.25}} & 8.23 \\
Benchmark suitable & 33 & 98.09 & 5.78 & \underline{\textbf{74.58}} & \underline{\textbf{4.71}} \\
\bottomrule
\end{tabular}
\end{table}

\textbf{Discussion. } Despite the fact that random outperforms SLNS on most tags, our approach can be advantageous in a specific context. The ML model gets a slight edge over random when the number of solutions found during probing is relatively high. That is compatible with the intuition that the prediction accuracy tends to improve with the quality of data. When the sample size is high, the features contain more signal that LGBM can take advantage of. Compared to SCIP, the probing data from Gurobi appears to be more informative.  All things considered, results are particularly noisy, which prevents us from making conclusive statements or strong recommendations. 

\subsection{LAP Results} \label{sec:lns:lap_results}

For the final set of experiments of this paper, we study the performance of the LNS matheuristic on the LAP. We compare against the SOA heuristic for this problem, which is PNF~\citep{robinson_la_rocca_one-shot_2024}.  The results in Table~\ref{tab:lns:lap_primal_gap} show the performance on 10 large LAP instances. Experiments have a 6-hour total time budget, a 30-minute budget for probing, a 15-minute budget per iteration and they use Gurobi for the MIP solver.  On top of the PNF results, we add the performance of random, OLNS and DOLNS. Oracles use the best solution found by Gurobi after 6 hours. Note that the new results use Algorithm~\ref{alg:lns}, which iteratively builds neighbourhoods instead of building a single one. The best algorithm is PNF-0.5-30 with a geometric mean of 1.07 and 5 wins. It is followed by PNFT-0.1-30, a threshold-based variant of PNF, with a geometric mean of 1.14 and 3 wins. The deterministic oracle DOLNS is the best performing algorithm that uses the LNS matheuristic with a geomean of 2.0. DOLNS outperforms OLNS because we did not add any noise to the oracles. We do not experiment with SLNS because the oracles fail to improve upon PNF. 

\begin{table}[h]
\centering
\caption{Primal gap for large LAP instances with 6-hour time limit}
\label{tab:lns:lap_primal_gap}
\begin{tabular}{lccccccc}
\toprule
{Scenario} & \multicolumn{3}{c}{Quantiles} & {Mean} & {Geomean} & {Wins} \\
{} & {0.1} & {0.5} & {0.9} & {} & {} & {} \\
\midrule
random & 1.48 & 1.81 & 3.33 & 2.07 & 2.96 & 0 \\
OLNS & 0.21 & 1.10 & 11.90 & 4.63 & 3.32 & 0 \\
DOLNS & 0.24 & 0.49 & 5.87 & 1.77 & 2.00 & 0 \\
grb-f1 & 0.46 & 0.97 & 7.25 & 2.44 & 2.53 & 0 \\
PNFT-0.1-30 & 0.00 & 0.10 & 0.33 & 0.15 & 1.14 & 3 \\
PNF-0.2-30 & 0.00 & 4.73 & 8.97 & 4.56 & 3.76 & 1 \\
PNF-0.5-30 & 0.00 & 0.01 & 0.24 & 0.08 & 1.07 & 5 \\
RINS-0.2-30 & 1.34 & 1.51 & 5.62 & 2.88 & 3.14 & 0 \\
RINS-0.5-30 & 4.72 & 5.08 & 5.45 & 5.09 & 6.08 & 0 \\
LB-0.7-30 & 0.11 & 10.15 & 10.92 & 7.41 & 5.80 & 1 \\
grb & 4.05 & 4.98 & 5.38 & 4.83 & 5.80 & 0 \\
\bottomrule
\end{tabular}
\end{table}

\textbf{Discussion. } 
There are at least three reasons that explain why oracles do not outperform PNF. First, the label solution found after 6 hours is not optimal on these large instances. Second, there is a 15-minute time limit on each LNS iteration, and this is not always enough for the solver to recover the label solution. Third, PNF fixes all the variables in the selected SOS1 constraints, whereas the oracles fix variables randomly across all constraints.  In other words, PNF explicitly exploit the structure of the assignments in LAP to produce a more natural sub-problem that is easier to solve.  Nevertheless, the LNS matheuristic is competitive for how general it is. For instance, random LNS is the third-best algorithm in terms of the 90$^{\text{th}}$ quantile (3.33) and the average primal gap (2.07). 


\section{Conclusion} \label{sec:lns:conclusion}

In this paper, our goal was to study how supervised learning can be used to accelerate the solution process of MIPs. To that end, we implemented an LNS matheuristic with a modular and flexible design. This allowed us to experiment with different LNS heuristics that produce neighbourhoods defined on the set of binary variables in the given instance. To test if our method is effective in general, we centered the experiments around MIPLIB, a well-known dataset of real-world instances from various domains. 

In our methodology, we proposed to represent the variable selection task as a supervised learning problem. For each binary variable, we attempt to predict the corresponding label solution with the help of a binary classifier. To validate this methodology, we experimented with different configurations of oracles that simulate the behaviour of the noisy ML predictions. We discovered a clear correlation between the error rate and the quality of solutions found at inference. We also found that the stochastic variant is more resilient to noise. This is a noteworthy insight because it implies one can improve test performance on MIPs with better prediction accuracy. 

We started our research effort with the goal to keep the training procedure minimal. To help improve the effectiveness of the ML model, we tested weight adjustments in the binary loss and two data sources for the features. Our results clearly show that loss adjustments can lower the FNR and they bring a noticeable performance improvement. The amount of data generated by SCIP was low compared to Gurobi, so we developed a sampling procedure for MIPs. The additional data from it was useful to compensate SCIP's shortcomings.  

The standout observation from this study is the strong performance of the random LNS. On different tags of MIPLIB and LAP, random LNS is competitive which is surprising given how simple it is. Prior research has not shown statistically significant improvements over random on the MIPLIB benchmark \citep{sonnerat_learning_2022, liu_machine_2023}. This suggests that developing a variable selection policy that consistently outperforms on a diverse set of unrelated MIPs remains a challenging endeavour. Two LNS baselines that we tested direct the search towards regions near the continuous relaxation solution, which contain sub-optimal solutions most of the time. Our heuristic, SLNS, outperformed all other non-random baselines but more research is required to get around its limitations. It lacks the ability to adapt to the instance on which it is tested. For example, one specific level of bias adjustment cannot be appropriate for all instances, since the ratio of 0 and 1 in the label solution can vary significantly from one instance to the other. Furthermore, the sampling routine is rather naive and it does not necessarily capture useful information for all instances. In order to realize substantial performance gains, we must improve the accuracy of the ML model. This sets up the path for further research.

\textbf{Future work. } One possible path for future work is to scale to large models~\citep{kaplan_scaling_2020}, which is considered a reliable way to improve the performance of ML models. We suspect, however, that scaling the model alone will not be enough to get appreciable improvements with the current dataset.  MIPLIB contains in total 1,065 instances but the majority can be solved to near-optimality within 10 minutes with the SOA commercial solver. We are left with, at most, a few hundred instances if we consider only the difficult ones. SOA models in image recognition and natural language processing are trained on datasets that are many orders of magnitude larger. For example, the image model SEER~\citep{goyal_self-supervised_2021} was trained in a self-supervised way on 1 billion images and the language model Llama 3~\citep{dubey_llama_2024_fixed} was trained on approximately 15.3 trillion tokens. To produce a large dataset, the established strategy is to use data augmentation techniques that create synthetic examples. Cropping and flipping are examples of image data augmentation techniques that modify the original image without changing its label~\citep{shorten_survey_2019}. For MIPs, data augmentation is more nuanced because small changes to its original structure can have a large impact on its properties. How can we design data augmentation techniques for MIPs that preserve the inherent characteristics of real-world problems? How can we ensure that synthetically generated MIP datasets accurately replicate the distribution of original instances? These are examples of questions that require more research effort to be answered.  We believe that a major milestone necessary to obtain appreciable performance improvements on generic MIPs is the generation of a large dataset of representative instances. 


\newpage

\clearpage

\newpage

\newpage

\clearpage

\newpage

\bibliographystyle{elsarticle-harv}

\bibliography{all_bib}

\begin{thebibliography}{44}
\expandafter\ifx\csname natexlab\endcsname\relax\def\natexlab#1{#1}\fi
\providecommand{\url}[1]{\texttt{#1}}
\providecommand{\href}[2]{#2}
\providecommand{\path}[1]{#1}
\providecommand{\DOIprefix}{doi:}
\providecommand{\ArXivprefix}{arXiv:}
\providecommand{\URLprefix}{URL: }
\providecommand{\Pubmedprefix}{pmid:}
\providecommand{\doi}[1]{\href{http://dx.doi.org/#1}{\path{#1}}}
\providecommand{\Pubmed}[1]{\href{pmid:#1}{\path{#1}}}
\providecommand{\bibinfo}[2]{#2}
\ifx\xfnm\relax \def\xfnm[#1]{\unskip,\space#1}\fi
\bibitem[{Bengio et~al.(2021)Bengio, Lodi and Prouvost}]{bengio_machine_2021}
\bibinfo{author}{Bengio, Y.}, \bibinfo{author}{Lodi, A.},
  \bibinfo{author}{Prouvost, A.}, \bibinfo{year}{2021}.
\newblock \bibinfo{title}{Machine learning for combinatorial optimization: {A}
  methodological tour d’horizon}.
\newblock \bibinfo{journal}{European Journal of Operational Research}
  \bibinfo{volume}{290}, \bibinfo{pages}{405--421}.
\newblock \DOIprefix\doi{10.1016/j.ejor.2020.07.063}.
\bibitem[{Bertacco et~al.(2007)Bertacco, Fischetti and
  Lodi}]{bertacco_feasibility_2007}
\bibinfo{author}{Bertacco, L.}, \bibinfo{author}{Fischetti, M.},
  \bibinfo{author}{Lodi, A.}, \bibinfo{year}{2007}.
\newblock \bibinfo{title}{A feasibility pump heuristic for general
  mixed-integer problems}.
\newblock \bibinfo{journal}{Discrete Optimization} \bibinfo{volume}{4},
  \bibinfo{pages}{63--76}.
\newblock \DOIprefix\doi{10.1016/j.disopt.2006.10.001}.
\bibitem[{Bestuzheva et~al.(2021)Bestuzheva, Besançon, Chen, Chmiela,
  Donkiewicz, van Doornmalen, Eifler, Gaul, Gamrath, Gleixner, Gottwald,
  Graczyk, Halbig, Hoen, Hojny, van~der Hulst, Koch, Lübbecke, Maher, Matter,
  Mühmer, Müller, Pfetsch, Rehfeldt, Schlein, Schl\"{o}sser, Serrano,
  Shinano, Sofranac, Turner, Vigerske, Wegscheider, Wellner, Weninger and
  Witzig}]{bestuzheva_scip_2021}
\bibinfo{author}{Bestuzheva, K.}, \bibinfo{author}{Besançon, M.},
  \bibinfo{author}{Chen, W.K.}, \bibinfo{author}{Chmiela, A.},
  \bibinfo{author}{Donkiewicz, T.}, \bibinfo{author}{van Doornmalen, J.},
  \bibinfo{author}{Eifler, L.}, \bibinfo{author}{Gaul, O.},
  \bibinfo{author}{Gamrath, G.}, \bibinfo{author}{Gleixner, A.},
  \bibinfo{author}{Gottwald, L.}, \bibinfo{author}{Graczyk, C.},
  \bibinfo{author}{Halbig, K.}, \bibinfo{author}{Hoen, A.},
  \bibinfo{author}{Hojny, C.}, \bibinfo{author}{van~der Hulst, R.},
  \bibinfo{author}{Koch, T.}, \bibinfo{author}{Lübbecke, M.},
  \bibinfo{author}{Maher, S.J.}, \bibinfo{author}{Matter, F.},
  \bibinfo{author}{Mühmer, E.}, \bibinfo{author}{Müller, B.},
  \bibinfo{author}{Pfetsch, M.E.}, \bibinfo{author}{Rehfeldt, D.},
  \bibinfo{author}{Schlein, S.}, \bibinfo{author}{Schl\"{o}sser, F.},
  \bibinfo{author}{Serrano, F.}, \bibinfo{author}{Shinano, Y.},
  \bibinfo{author}{Sofranac, B.}, \bibinfo{author}{Turner, M.},
  \bibinfo{author}{Vigerske, S.}, \bibinfo{author}{Wegscheider, F.},
  \bibinfo{author}{Wellner, P.}, \bibinfo{author}{Weninger, D.},
  \bibinfo{author}{Witzig, J.}, \bibinfo{year}{2021}.
\newblock \bibinfo{title}{The {SCIP} {Optimization} {Suite} 8.0}.
\newblock \DOIprefix\doi{10.48550/arXiv.2112.08872}. \bibinfo{note}{arXiv
  preprint arXiv:2112.08872}.
\bibitem[{Bolusani et~al.(2024)Bolusani, Besançon, Bestuzheva, Chmiela,
  Dionísio, Donkiewicz, van Doornmalen, Eifler, Ghannam, Gleixner, Graczyk,
  Halbig, Hedtke, Hoen, Hojny, van~der Hulst, Kamp, Koch, Kofler, Lentz, Manns,
  Mexi, Mühmer, Pfetsch, Schl\"{o}sser, Serrano, Shinano, Turner, Vigerske,
  Weninger and Xu}]{bolusani_scip_2024}
\bibinfo{author}{Bolusani, S.}, \bibinfo{author}{Besançon, M.},
  \bibinfo{author}{Bestuzheva, K.}, \bibinfo{author}{Chmiela, A.},
  \bibinfo{author}{Dionísio, J.}, \bibinfo{author}{Donkiewicz, T.},
  \bibinfo{author}{van Doornmalen, J.}, \bibinfo{author}{Eifler, L.},
  \bibinfo{author}{Ghannam, M.}, \bibinfo{author}{Gleixner, A.},
  \bibinfo{author}{Graczyk, C.}, \bibinfo{author}{Halbig, K.},
  \bibinfo{author}{Hedtke, I.}, \bibinfo{author}{Hoen, A.},
  \bibinfo{author}{Hojny, C.}, \bibinfo{author}{van~der Hulst, R.},
  \bibinfo{author}{Kamp, D.}, \bibinfo{author}{Koch, T.},
  \bibinfo{author}{Kofler, K.}, \bibinfo{author}{Lentz, J.},
  \bibinfo{author}{Manns, J.}, \bibinfo{author}{Mexi, G.},
  \bibinfo{author}{Mühmer, E.}, \bibinfo{author}{Pfetsch, M.E.},
  \bibinfo{author}{Schl\"{o}sser, F.}, \bibinfo{author}{Serrano, F.},
  \bibinfo{author}{Shinano, Y.}, \bibinfo{author}{Turner, M.},
  \bibinfo{author}{Vigerske, S.}, \bibinfo{author}{Weninger, D.},
  \bibinfo{author}{Xu, L.}, \bibinfo{year}{2024}.
\newblock \bibinfo{title}{The {SCIP} {Optimization} {Suite} 9.0}.
\newblock \DOIprefix\doi{10.48550/arXiv.2402.17702}. \bibinfo{note}{arXiv
  preprint arXiv:2402.17702}.
\bibitem[{Borisov et~al.(2022)Borisov, Leemann, Seßler, Haug, Pawelczyk and
  Kasneci}]{borisov_deep_2022_2}
\bibinfo{author}{Borisov, V.}, \bibinfo{author}{Leemann, T.},
  \bibinfo{author}{Seßler, K.}, \bibinfo{author}{Haug, J.},
  \bibinfo{author}{Pawelczyk, M.}, \bibinfo{author}{Kasneci, G.},
  \bibinfo{year}{2022}.
\newblock \bibinfo{title}{Deep {Neural} {Networks} and {Tabular} {Data}: {A}
  {Survey}}.
\newblock \bibinfo{journal}{IEEE Transactions on Neural Networks and Learning
  Systems} \bibinfo{volume}{1}, \bibinfo{pages}{1--21}.
\newblock \DOIprefix\doi{10.1109/TNNLS.2022.3229161}.
\bibitem[{Caruana and Niculescu-Mizil(2006)}]{caruana_empirical_2006}
\bibinfo{author}{Caruana, R.}, \bibinfo{author}{Niculescu-Mizil, A.},
  \bibinfo{year}{2006}.
\newblock \bibinfo{title}{An empirical comparison of supervised learning
  algorithms}, in: \bibinfo{booktitle}{Proceedings of the 23rd international
  conference on {Machine} learning}, \bibinfo{publisher}{Association for
  Computing Machinery}, \bibinfo{address}{New York, NY, USA}. pp.
  \bibinfo{pages}{161--168}.
\newblock \DOIprefix\doi{10.1145/1143844.1143865}.
\bibitem[{Chmiela et~al.(2023)Chmiela, Gleixner, Lichocki and
  Pokutta}]{chmiela_online_2023}
\bibinfo{author}{Chmiela, A.}, \bibinfo{author}{Gleixner, A.},
  \bibinfo{author}{Lichocki, P.}, \bibinfo{author}{Pokutta, S.},
  \bibinfo{year}{2023}.
\newblock \bibinfo{title}{Online {Learning} for {Scheduling} {MIP}
  {Heuristics}}, in: \bibinfo{editor}{Cire, A.A.} (Ed.),
  \bibinfo{booktitle}{Integration of {Constraint} {Programming}, {Artificial}
  {Intelligence}, and {Operations} {Research}}, \bibinfo{publisher}{Springer
  Nature Switzerland}, \bibinfo{address}{Cham}. pp. \bibinfo{pages}{114--123}.
\newblock \DOIprefix\doi{10.1007/978-3-031-33271-5-8}.
\bibitem[{Crainic et~al.(2004)Crainic, Gendron and Hernu}]{crainic_slope_2004}
\bibinfo{author}{Crainic, T.G.}, \bibinfo{author}{Gendron, B.},
  \bibinfo{author}{Hernu, G.}, \bibinfo{year}{2004}.
\newblock \bibinfo{title}{A {Slope} {Scaling}/{Lagrangean} {Perturbation}
  {Heuristic} with {Long}-{Term} {Memory} for {Multicommodity} {Capacitated}
  {Fixed}-{Charge} {Network} {Design}}.
\newblock \bibinfo{journal}{Journal of Heuristics} \bibinfo{volume}{10},
  \bibinfo{pages}{525--545}.
\newblock \DOIprefix\doi{10.1023/B:HEUR.0000045323.83583.bd}.
\bibitem[{Danna et~al.(2005)Danna, Rothberg and Pape}]{danna_exploring_2005}
\bibinfo{author}{Danna, E.}, \bibinfo{author}{Rothberg, E.},
  \bibinfo{author}{Pape, C.L.}, \bibinfo{year}{2005}.
\newblock \bibinfo{title}{Exploring relaxation induced neighborhoods to improve
  {MIP} solutions}.
\newblock \bibinfo{journal}{Mathematical Programming} \bibinfo{volume}{102},
  \bibinfo{pages}{71--90}.
\newblock \DOIprefix\doi{10.1007/s10107-004-0518-7}.
\bibitem[{Dantzig(1951)}]{dantzig1951maximization}
\bibinfo{author}{Dantzig, G.B.}, \bibinfo{year}{1951}.
\newblock \bibinfo{title}{Maximization of a linear function of variables
  subject to linear inequalities}.
\newblock \bibinfo{journal}{Activity analysis of production and allocation}
  \bibinfo{volume}{13}, \bibinfo{pages}{339--347}.
\bibitem[{Dubey et~al.(2024)Dubey, Jauhri and et~al.}]{dubey_llama_2024_fixed}
\bibinfo{author}{Dubey, A.}, \bibinfo{author}{Jauhri, A.},
  \bibinfo{author}{et~al.}, \bibinfo{year}{2024}.
\newblock \bibinfo{title}{The {Llama} 3 {Herd} of {Models}}.
\newblock \DOIprefix\doi{10.48550/arXiv.2407.21783}. \bibinfo{note}{arXiv
  preprint arXiv:2407.21783}.
\bibitem[{Feijen et~al.(2024)Feijen, Schäfer, Dekker and
  Pieterse}]{feijen_learning-enhanced_2024}
\bibinfo{author}{Feijen, W.}, \bibinfo{author}{Schäfer, G.},
  \bibinfo{author}{Dekker, K.}, \bibinfo{author}{Pieterse, S.},
  \bibinfo{year}{2024}.
\newblock \bibinfo{title}{Learning-{Enhanced} {Neighborhood} {Selection} for
  the {Vehicle} {Routing} {Problem} with {Time} {Windows}}.
\newblock \DOIprefix\doi{10.48550/arXiv.2403.08839}. \bibinfo{note}{arXiv
  preprint arXiv:2403.08839}.
\bibitem[{Fischetti and Fischetti(2018)}]{fischetti_matheuristics_2018}
\bibinfo{author}{Fischetti, M.}, \bibinfo{author}{Fischetti, M.},
  \bibinfo{year}{2018}.
\newblock \bibinfo{title}{Matheuristics}, in: \bibinfo{editor}{Martí, R.},
  \bibinfo{editor}{Pardalos, P.M.}, \bibinfo{editor}{Resende, M.G.C.} (Eds.),
  \bibinfo{booktitle}{Handbook of {Heuristics}}. \bibinfo{publisher}{Springer
  International Publishing}, \bibinfo{address}{Cham}, pp.
  \bibinfo{pages}{121--153}.
\newblock \DOIprefix\doi{10.1007/978-3-319-07124-4-14}.
\bibitem[{Fischetti et~al.(2005)Fischetti, Glover and
  Lodi}]{fischetti_feasibility_2005}
\bibinfo{author}{Fischetti, M.}, \bibinfo{author}{Glover, F.},
  \bibinfo{author}{Lodi, A.}, \bibinfo{year}{2005}.
\newblock \bibinfo{title}{The feasibility pump}.
\newblock \bibinfo{journal}{Mathematical Programming} \bibinfo{volume}{104},
  \bibinfo{pages}{91--104}.
\newblock \DOIprefix\doi{10.1007/s10107-004-0570-3}.
\bibitem[{Gleixner et~al.(2021)Gleixner, Hendel, Gamrath, Achterberg, Bastubbe,
  Berthold, Christophel, Jarck, Koch, Linderoth, Lübbecke, Mittelmann, Ozyurt,
  Ralphs, Salvagnin and Shinano}]{gleixner_miplib_2021}
\bibinfo{author}{Gleixner, A.}, \bibinfo{author}{Hendel, G.},
  \bibinfo{author}{Gamrath, G.}, \bibinfo{author}{Achterberg, T.},
  \bibinfo{author}{Bastubbe, M.}, \bibinfo{author}{Berthold, T.},
  \bibinfo{author}{Christophel, P.}, \bibinfo{author}{Jarck, K.},
  \bibinfo{author}{Koch, T.}, \bibinfo{author}{Linderoth, J.},
  \bibinfo{author}{Lübbecke, M.}, \bibinfo{author}{Mittelmann, H.D.},
  \bibinfo{author}{Ozyurt, D.}, \bibinfo{author}{Ralphs, T.K.},
  \bibinfo{author}{Salvagnin, D.}, \bibinfo{author}{Shinano, Y.},
  \bibinfo{year}{2021}.
\newblock \bibinfo{title}{{MIPLIB} 2017: data-driven compilation of the
  6th mixed-integer programming library}.
\newblock \bibinfo{journal}{Mathematical Programming Computation}
  \bibinfo{volume}{13}, \bibinfo{pages}{443--490}.
\newblock \DOIprefix\doi{10.1007/s12532-020-00194-3}.
\bibitem[{Goyal et~al.(2021)Goyal, Caron, Lefaudeux, Xu, Wang, Pai, Singh,
  Liptchinsky, Misra, Joulin and Bojanowski}]{goyal_self-supervised_2021}
\bibinfo{author}{Goyal, P.}, \bibinfo{author}{Caron, M.},
  \bibinfo{author}{Lefaudeux, B.}, \bibinfo{author}{Xu, M.},
  \bibinfo{author}{Wang, P.}, \bibinfo{author}{Pai, V.},
  \bibinfo{author}{Singh, M.}, \bibinfo{author}{Liptchinsky, V.},
  \bibinfo{author}{Misra, I.}, \bibinfo{author}{Joulin, A.},
  \bibinfo{author}{Bojanowski, P.}, \bibinfo{year}{2021}.
\newblock \bibinfo{title}{Self-supervised {Pretraining} of {Visual} {Features}
  in the {Wild}}.
\newblock \DOIprefix\doi{10.48550/arXiv.2103.01988}. \bibinfo{note}{arXiv
  preprint arXiv:2103.01988}.
\bibitem[{{Gurobi Optimization, LLC}(2024)}]{gurobi}
\bibinfo{author}{{Gurobi Optimization, LLC}}, \bibinfo{year}{2024}.
\newblock \bibinfo{title}{{Gurobi Optimizer Reference Manual}}.
\newblock \URLprefix \url{https://www.gurobi.com}. \bibinfo{note}{{Accessed:
  2024-10-23}}.
\bibitem[{Hendel(2022)}]{hendel_adaptive_2022}
\bibinfo{author}{Hendel, G.}, \bibinfo{year}{2022}.
\newblock \bibinfo{title}{Adaptive large neighborhood search for mixed integer
  programming}.
\newblock \bibinfo{journal}{Mathematical Programming Computation}
  \bibinfo{volume}{14}, \bibinfo{pages}{185--221}.
\newblock \DOIprefix\doi{10.1007/s12532-021-00209-7}.
\bibitem[{Hendel et~al.(2019)Hendel, Miltenberger and
  Witzig}]{hendel_adaptive_2019}
\bibinfo{author}{Hendel, G.}, \bibinfo{author}{Miltenberger, M.},
  \bibinfo{author}{Witzig, J.}, \bibinfo{year}{2019}.
\newblock \bibinfo{title}{Adaptive {Algorithmic} {Behavior} for {Solving}
  {Mixed} {Integer} {Programs} {Using} {Bandit} {Algorithms}}, in:
  \bibinfo{editor}{Fortz, B.}, \bibinfo{editor}{Labbé, M.} (Eds.),
  \bibinfo{booktitle}{Operations {Research} {Proceedings} 2018},
  \bibinfo{publisher}{Springer International Publishing},
  \bibinfo{address}{Cham}. pp. \bibinfo{pages}{513--519}.
\newblock \DOIprefix\doi{10.1007/978-3-030-18500-8-64}.
\bibitem[{Huang et~al.(2023a)Huang, Ferber, Tian, Dilkina and
  Steiner}]{huang_local_2023}
\bibinfo{author}{Huang, T.}, \bibinfo{author}{Ferber, A.},
  \bibinfo{author}{Tian, Y.}, \bibinfo{author}{Dilkina, B.},
  \bibinfo{author}{Steiner, B.}, \bibinfo{year}{2023}a.
\newblock \bibinfo{title}{Local {Branching} {Relaxation} {Heuristics}
  for {Integer} {Linear} {Programs}}, in: \bibinfo{editor}{Cire, A.A.} (Ed.),
  \bibinfo{booktitle}{Integration of {Constraint} {Programming}, {Artificial}
  {Intelligence}, and {Operations} {Research}}, \bibinfo{publisher}{Springer
  Nature Switzerland}, \bibinfo{address}{Cham}. pp. \bibinfo{pages}{96--113}.
\newblock \DOIprefix\doi{10.1007/978-3-031-33271-5-7}.
\bibitem[{Huang et~al.(2023b)Huang, Ferber, Tian, Dilkina and
  Steiner}]{huang_searching_2023}
\bibinfo{author}{Huang, T.}, \bibinfo{author}{Ferber, A.M.},
  \bibinfo{author}{Tian, Y.}, \bibinfo{author}{Dilkina, B.},
  \bibinfo{author}{Steiner, B.}, \bibinfo{year}{2023}b.
\newblock \bibinfo{title}{Searching {Large} {Neighborhoods} for {Integer}
  {Linear} {Programs} with {Contrastive} {Learning}}, in:
  \bibinfo{booktitle}{Proceedings of the 40th {International} {Conference} on
  {Machine} {Learning}}, \bibinfo{publisher}{PMLR}. pp.
  \bibinfo{pages}{13869--13890}.
\bibitem[{Huangfu and Hall(2018)}]{huangfu_parallelizing_2018}
\bibinfo{author}{Huangfu, Q.}, \bibinfo{author}{Hall, J.A.J.},
  \bibinfo{year}{2018}.
\newblock \bibinfo{title}{Parallelizing the dual revised simplex method}.
\newblock \bibinfo{journal}{Mathematical Programming Computation}
  \bibinfo{volume}{10}, \bibinfo{pages}{119--142}.
\newblock \DOIprefix\doi{10.1007/s12532-017-0130-5}.
\bibitem[{Kaplan et~al.(2020)Kaplan, McCandlish, Henighan, Brown, Chess, Child,
  Gray, Radford, Wu and Amodei}]{kaplan_scaling_2020}
\bibinfo{author}{Kaplan, J.}, \bibinfo{author}{McCandlish, S.},
  \bibinfo{author}{Henighan, T.}, \bibinfo{author}{Brown, T.B.},
  \bibinfo{author}{Chess, B.}, \bibinfo{author}{Child, R.},
  \bibinfo{author}{Gray, S.}, \bibinfo{author}{Radford, A.},
  \bibinfo{author}{Wu, J.}, \bibinfo{author}{Amodei, D.}, \bibinfo{year}{2020}.
\newblock \bibinfo{title}{Scaling {Laws} for {Neural} {Language} {Models}}.
\newblock \DOIprefix\doi{10.48550/arXiv.2001.08361}. \bibinfo{note}{arXiv
  preprint arXiv:2001.08361}.
\bibitem[{Karmarkar(1984)}]{karmarkar_new_1984}
\bibinfo{author}{Karmarkar, N.}, \bibinfo{year}{1984}.
\newblock \bibinfo{title}{A new polynomial-time algorithm for linear
  programming}, in: \bibinfo{booktitle}{Proceedings of the sixteenth annual
  {ACM} symposium on {Theory} of computing}, \bibinfo{publisher}{Association
  for Computing Machinery}, \bibinfo{address}{New York, NY, USA}. pp.
  \bibinfo{pages}{302--311}.
\newblock \DOIprefix\doi{10.1145/800057.808695}.
\bibitem[{Ke et~al.(2017)Ke, Meng, Finley, Wang, Chen, Ma, Ye and
  Liu}]{ke_lightgbm_2017}
\bibinfo{author}{Ke, G.}, \bibinfo{author}{Meng, Q.}, \bibinfo{author}{Finley,
  T.}, \bibinfo{author}{Wang, T.}, \bibinfo{author}{Chen, W.},
  \bibinfo{author}{Ma, W.}, \bibinfo{author}{Ye, Q.}, \bibinfo{author}{Liu,
  T.Y.}, \bibinfo{year}{2017}.
\newblock \bibinfo{title}{{LightGBM}: {A} {Highly} {Efficient} {Gradient}
  {Boosting} {Decision} {Tree}}, in: \bibinfo{booktitle}{Advances in {Neural}
  {Information} {Processing} {Systems}}, \bibinfo{publisher}{Curran Associates,
  Inc.}
\bibitem[{Khalil et~al.(2016)Khalil, Bodic, Song, Nemhauser and
  Dilkina}]{khalil_learning_2016}
\bibinfo{author}{Khalil, E.}, \bibinfo{author}{Bodic, P.L.},
  \bibinfo{author}{Song, L.}, \bibinfo{author}{Nemhauser, G.},
  \bibinfo{author}{Dilkina, B.}, \bibinfo{year}{2016}.
\newblock \bibinfo{title}{Learning to {Branch} in {Mixed} {Integer}
  {Programming}}.
\newblock \bibinfo{journal}{Proceedings of the AAAI Conference on Artificial
  Intelligence} \bibinfo{volume}{30}.
\newblock \DOIprefix\doi{10.1609/aaai.v30i1.10080}.
\bibitem[{Khalil et~al.(2022)Khalil, Morris and Lodi}]{khalil_mip-gnn_2022}
\bibinfo{author}{Khalil, E.B.}, \bibinfo{author}{Morris, C.},
  \bibinfo{author}{Lodi, A.}, \bibinfo{year}{2022}.
\newblock \bibinfo{title}{{MIP}-{GNN}: {A} {Data}-{Driven} {Framework} for
  {Guiding} {Combinatorial} {Solvers}}.
\newblock \bibinfo{journal}{Proceedings of the AAAI Conference on Artificial
  Intelligence} \bibinfo{volume}{36}, \bibinfo{pages}{10219--10227}.
\newblock \DOIprefix\doi{10.1609/aaai.v36i9.21262}.
\bibitem[{Kirkpatrick et~al.(1983)Kirkpatrick, Gelatt and
  Vecchi}]{kirkpatrick_optimization_1983}
\bibinfo{author}{Kirkpatrick, S.}, \bibinfo{author}{Gelatt, C.D.},
  \bibinfo{author}{Vecchi, M.P.}, \bibinfo{year}{1983}.
\newblock \bibinfo{title}{Optimization by {Simulated} {Annealing}}.
\newblock \bibinfo{journal}{Science} \bibinfo{volume}{220},
  \bibinfo{pages}{671--680}.
\newblock \DOIprefix\doi{10.1126/science.220.4598.671}.
\bibitem[{Kotary et~al.(2021)Kotary, Fioretto and
  Van~Hentenryck}]{kotary_learning_2021}
\bibinfo{author}{Kotary, J.}, \bibinfo{author}{Fioretto, F.},
  \bibinfo{author}{Van~Hentenryck, P.}, \bibinfo{year}{2021}.
\newblock \bibinfo{title}{Learning hard optimization problems: {A} data
  generation perspective}.
\newblock \bibinfo{journal}{Advances in Neural Information Processing Systems}
  \bibinfo{volume}{34}, \bibinfo{pages}{24981--24992}.
\bibitem[{Liu et~al.(2022)Liu, Fischetti and Lodi}]{liu_learning_2022}
\bibinfo{author}{Liu, D.}, \bibinfo{author}{Fischetti, M.},
  \bibinfo{author}{Lodi, A.}, \bibinfo{year}{2022}.
\newblock \bibinfo{title}{Learning to {Search} in {Local} {Branching}}.
\newblock \bibinfo{journal}{Proceedings of the AAAI Conference on Artificial
  Intelligence} \bibinfo{volume}{36}, \bibinfo{pages}{3796--3803}.
\newblock \DOIprefix\doi{10.1609/aaai.v36i4.20294}.
\bibitem[{Liu et~al.(2023)Liu, Perreault, Hertz and Lodi}]{liu_machine_2023}
\bibinfo{author}{Liu, D.}, \bibinfo{author}{Perreault, V.},
  \bibinfo{author}{Hertz, A.}, \bibinfo{author}{Lodi, A.},
  \bibinfo{year}{2023}.
\newblock \bibinfo{title}{A machine learning framework for neighbor generation
  in metaheuristic search}.
\newblock \bibinfo{journal}{Frontiers in Applied Mathematics and Statistics}
  \bibinfo{volume}{9}.
\newblock \DOIprefix\doi{10.3389/fams.2023.1128181}.
\bibitem[{Lubin et~al.(2023)Lubin, Dowson, Garcia, Huchette, Legat and
  Vielma}]{lubin_jump_2023}
\bibinfo{author}{Lubin, M.}, \bibinfo{author}{Dowson, O.},
  \bibinfo{author}{Garcia, J.D.}, \bibinfo{author}{Huchette, J.},
  \bibinfo{author}{Legat, B.}, \bibinfo{author}{Vielma, J.P.},
  \bibinfo{year}{2023}.
\newblock \bibinfo{title}{{JuMP} 1.0: recent improvements to a modeling
  language for mathematical optimization}.
\newblock \bibinfo{journal}{Mathematical Programming Computation}
  \bibinfo{volume}{15}, \bibinfo{pages}{581--589}.
\newblock \DOIprefix\doi{10.1007/s12532-023-00239-3}.
\bibitem[{Nair et~al.(2021)Nair, Bartunov, Gimeno, von Glehn, Lichocki, Lobov,
  O'Donoghue, Sonnerat, Tjandraatmadja, Wang, Addanki, Hapuarachchi, Keck,
  Keeling, Kohli, Ktena, Li, Vinyals and Zwols}]{nair_solving_2021}
\bibinfo{author}{Nair, V.}, \bibinfo{author}{Bartunov, S.},
  \bibinfo{author}{Gimeno, F.}, \bibinfo{author}{von Glehn, I.},
  \bibinfo{author}{Lichocki, P.}, \bibinfo{author}{Lobov, I.},
  \bibinfo{author}{O'Donoghue, B.}, \bibinfo{author}{Sonnerat, N.},
  \bibinfo{author}{Tjandraatmadja, C.}, \bibinfo{author}{Wang, P.},
  \bibinfo{author}{Addanki, R.}, \bibinfo{author}{Hapuarachchi, T.},
  \bibinfo{author}{Keck, T.}, \bibinfo{author}{Keeling, J.},
  \bibinfo{author}{Kohli, P.}, \bibinfo{author}{Ktena, I.},
  \bibinfo{author}{Li, Y.}, \bibinfo{author}{Vinyals, O.},
  \bibinfo{author}{Zwols, Y.}, \bibinfo{year}{2021}.
\newblock \bibinfo{title}{Solving {Mixed} {Integer} {Programs} {Using} {Neural}
  {Networks}}.
\newblock \DOIprefix\doi{10.48550/arXiv.2012.13349}. \bibinfo{note}{arXiv
  preprint arXiv:2012.13349}.
\bibitem[{Ortiz-Astorquiza et~al.(2021)Ortiz-Astorquiza, Cordeau and
  Frejinger}]{ortiz-astorquiza_locomotive_2021}
\bibinfo{author}{Ortiz-Astorquiza, C.}, \bibinfo{author}{Cordeau, {\relax
  J.-F}.}, \bibinfo{author}{Frejinger, E.}, \bibinfo{year}{2021}.
\newblock \bibinfo{title}{The {Locomotive} {Assignment} {Problem} with
  {Distributed} {Power} at the {Canadian} {National} {Railway} {Company}}.
\newblock \bibinfo{journal}{Transportation Science} \bibinfo{volume}{55},
  \bibinfo{pages}{510--531}.
\newblock \DOIprefix\doi{10.1287/trsc.2020.1030}.
\bibitem[{Peng et~al.(2021)Peng, Choi and Xu}]{peng_graph_2021}
\bibinfo{author}{Peng, Y.}, \bibinfo{author}{Choi, B.}, \bibinfo{author}{Xu,
  J.}, \bibinfo{year}{2021}.
\newblock \bibinfo{title}{Graph {Learning} for {Combinatorial} {Optimization}:
  {A} {Survey} of {State}-of-the-{Art}}.
\newblock \bibinfo{journal}{Data Science and Engineering} \bibinfo{volume}{6},
  \bibinfo{pages}{119--141}.
\newblock \DOIprefix\doi{10.1007/s41019-021-00155-3}.
\bibitem[{Robinson La~Rocca et~al.(2024a)Robinson La~Rocca, Cordeau and
  Frejinger}]{robinson_la_rocca_combining_2024}
\bibinfo{author}{Robinson La~Rocca, C.}, \bibinfo{author}{Cordeau, {\relax
  J.-F}.}, \bibinfo{author}{Frejinger, E.}, \bibinfo{year}{2024}a.
\newblock \bibinfo{title}{Combining supervised learning and local search for
  the multicommodity capacitated fixed-charge network design problem}.
\newblock \bibinfo{journal}{Transportation Research Part E: Logistics and
  Transportation Review} \bibinfo{volume}{192}, \bibinfo{pages}{103805}.
\newblock \DOIprefix\doi{10.1016/j.tre.2024.103805}.
\bibitem[{Robinson La~Rocca et~al.(2024b)Robinson La~Rocca, Cordeau and
  Frejinger}]{robinson_la_rocca_one-shot_2024}
\bibinfo{author}{Robinson La~Rocca, C.}, \bibinfo{author}{Cordeau, {\relax
  J.-F}.}, \bibinfo{author}{Frejinger, E.}, \bibinfo{year}{2024}b.
\newblock \bibinfo{title}{One-{Shot} {Learning} for {MIPs} with {SOS1}
  {Constraints}}.
\newblock \bibinfo{journal}{Operations Research Forum} \bibinfo{volume}{5},
  \bibinfo{pages}{57}.
\newblock \DOIprefix\doi{10.1007/s43069-024-00336-6}.
\bibitem[{Ropke and Pisinger(2006)}]{ropke_adaptive_2006}
\bibinfo{author}{Ropke, S.}, \bibinfo{author}{Pisinger, D.},
  \bibinfo{year}{2006}.
\newblock \bibinfo{title}{An {Adaptive} {Large} {Neighborhood} {Search}
  {Heuristic} for the {Pickup} and {Delivery} {Problem} with {Time} {Windows}}.
\newblock \bibinfo{journal}{Transportation Science} \bibinfo{volume}{40},
  \bibinfo{pages}{455--472}.
\newblock \DOIprefix\doi{10.1287/trsc.1050.0135}.
\bibitem[{Rothberg(2007)}]{rothberg_evolutionary_2007}
\bibinfo{author}{Rothberg, E.}, \bibinfo{year}{2007}.
\newblock \bibinfo{title}{An {Evolutionary} {Algorithm} for {Polishing} {Mixed}
  {Integer} {Programming} {Solutions}}.
\newblock \bibinfo{journal}{INFORMS Journal on Computing} \bibinfo{volume}{19},
  \bibinfo{pages}{534--541}.
\newblock \DOIprefix\doi{10.1287/ijoc.1060.0189}.
\bibitem[{Scavuzzo et~al.(2024)Scavuzzo, Aardal, Lodi and
  Yorke-Smith}]{scavuzzo_machine_2024}
\bibinfo{author}{Scavuzzo, L.}, \bibinfo{author}{Aardal, K.},
  \bibinfo{author}{Lodi, A.}, \bibinfo{author}{Yorke-Smith, N.},
  \bibinfo{year}{2024}.
\newblock \bibinfo{title}{Machine learning augmented branch and bound for mixed
  integer linear programming}.
\newblock \bibinfo{journal}{Mathematical Programming}
  \DOIprefix\doi{10.1007/s10107-024-02130-y}.
\bibitem[{Shaw(1998)}]{shaw_using_1998}
\bibinfo{author}{Shaw, P.}, \bibinfo{year}{1998}.
\newblock \bibinfo{title}{Using {Constraint} {Programming} and {Local} {Search}
  {Methods} to {Solve} {Vehicle} {Routing} {Problems}}, in:
  \bibinfo{editor}{Maher, M.}, \bibinfo{editor}{Puget, J.F.} (Eds.),
  \bibinfo{booktitle}{Principles and {Practice} of {Constraint} {Programming}
  — {CP98}}, \bibinfo{publisher}{Springer}, \bibinfo{address}{Berlin,
  Heidelberg}. pp. \bibinfo{pages}{417--431}.
\newblock \DOIprefix\doi{10.1007/3-540-49481-2-30}.
\bibitem[{Shorten and Khoshgoftaar(2019)}]{shorten_survey_2019}
\bibinfo{author}{Shorten, C.}, \bibinfo{author}{Khoshgoftaar, T.M.},
  \bibinfo{year}{2019}.
\newblock \bibinfo{title}{A survey on {Image} {Data} {Augmentation} for {Deep}
  {Learning}}.
\newblock \bibinfo{journal}{Journal of Big Data} \bibinfo{volume}{6},
  \bibinfo{pages}{60}.
\newblock \DOIprefix\doi{10.1186/s40537-019-0197-0}.
\bibitem[{Sonnerat et~al.(2022)Sonnerat, Wang, Ktena, Bartunov and
  Nair}]{sonnerat_learning_2022}
\bibinfo{author}{Sonnerat, N.}, \bibinfo{author}{Wang, P.},
  \bibinfo{author}{Ktena, I.}, \bibinfo{author}{Bartunov, S.},
  \bibinfo{author}{Nair, V.}, \bibinfo{year}{2022}.
\newblock \bibinfo{title}{Learning a {Large} {Neighborhood} {Search}
  {Algorithm} for {Mixed} {Integer} {Programs}}.
\newblock \DOIprefix\doi{10.48550/arXiv.2107.10201}. \bibinfo{note}{arXiv
  preprint arXiv:2107.10201}.
\bibitem[{Sun et~al.(2021)Sun, Ernst, Li and Weiner}]{sun_generalization_2021}
\bibinfo{author}{Sun, Y.}, \bibinfo{author}{Ernst, A.}, \bibinfo{author}{Li,
  X.}, \bibinfo{author}{Weiner, J.}, \bibinfo{year}{2021}.
\newblock \bibinfo{title}{Generalization of machine learning for problem
  reduction: {A} case study on travelling salesman problems}.
\newblock \bibinfo{journal}{OR Spectrum} \bibinfo{volume}{43},
  \bibinfo{pages}{607--633}.
\newblock \DOIprefix\doi{10.1007/s00291-020-00604-x}.

\end{thebibliography}

\end{document}